\documentclass[11pt,a4paper]{article}
\usepackage{ascmac,amsmath,amssymb,bm}
\usepackage{multirow,url,%theorem,
array}
\usepackage{amsthm}
\usepackage[dvipdfmx]{graphicx}
\usepackage{color}

\newtheorem{thm}{Theorem}[section]
\newtheorem{rem}[thm]{Remark}
\newtheorem{cor}[thm]{Corollary}
\newtheorem{lem}[thm]{Lemma}

\newtheorem{assum}[thm]{Assumption}

\usepackage[top=30truemm,bottom=30truemm,left=25truemm,right=25truemm]{geometry}

\newcommand{\argmin}{\operatornamewithlimits{argmin}}

\newcommand{\p}{\partial}

\newcommand{\mbbm}{\mathbb{M}}

\numberwithin{equation}{section}
\allowdisplaybreaks

\newcommand{\mca}{\mathcal{A}}

\newcommand{\mcc}{\mathcal{C}}

\newcommand{\mcf}{\mathcal{F}}

\newcommand{\mcl}{\mathcal{L}}

\newcommand{\mbbh}{\mathbb{H}}
\newcommand{\mbbn}{\mathbb{N}}
\newcommand{\mbbr}{\mathbb{R}}
\newcommand{\mbbrp}{\mathbb{R}_{+}}

\newcommand{\mbX}{\mathbf{X}}

\newcommand{\al}{\alpha}
\newcommand{\del}{\delta}

\newcommand{\sig}{\sigma}
\newcommand{\ep}{\epsilon}
\newcommand{\D}{\Delta}
\newcommand{\Sig}{\Sigma}
\newcommand{\gam}{\gamma}
\newcommand{\lam}{\lambda}
\newcommand{\Lam}{\Lambda}
\newcommand{\Gam}{\Gamma}

\newcommand{\cil}{\xrightarrow{\mcl}} % <- Convergence in law
\newcommand{\scl}{\xrightarrow{\mcl_{s}}} % <- Stable convergence in law
\newcommand{\cip}{\xrightarrow{p}} % <- Convergence in probability

\newcommand{\diag}{\mathrm{diag}}

\def\ds#1{\displaystyle{#1}}
\def\nn{\nonumber}

\def\sumj{\sum_{j=1}^{n}}

\def\E{\mathbb{E}}

\def\ulop#1{O^{\ast}_{p}(#1)}
\def\usop#1{o^{\ast}_{p}(#1)}

\def\tz{\theta_{0}}
\def\tes{\hat{\theta}_{n}}
\def\aes{\hat{\alpha}_{n}}
\def\bes{\hat{\beta}_{n}}
\def\xes{\hat{\xi}_{n}}

\title{Moment convergence in regularized estimation under multiple and mixed-rates asymptotics\footnote{This version: \today}}
\author{Hiroki Masuda\footnote{Faculty of Mathematics, Kyushu University. 744 Motooka Nishi-ku Fukuoka 819-0395, Japan.}
\and Yusuke Shimizu\footnote{Department of Mathematics, Josai University. 1-1 Keyakidai Sakado Saitama 350-0295, Japan. \newline \quad \quad Corresponding author:\texttt{yshimizu@josai.ac.jp}}}
\date{\empty}

\begin{document}
\setlength{\baselineskip}{4.5mm}

\maketitle

\begin{abstract}
In $M$-estimation under standard asymptotics, the weak convergence combined with the polynomial type large deviation estimate of the associated statistical random field Yoshida (2011) provides us with not only the asymptotic distribution of the associated $M$-estimator but also the convergence of its moments, the latter playing an important role in theoretical statistics.
In this paper, we study the above program for statistical random fields of multiple and also possibly mixed-rates type in the sense of Radchenko (2008) where the associated statistical random fields may be non-differentiable and may fail to be locally asymptotically quadratic. Consequently, a very strong mode of convergence of a wide range of regularized $M$-estimators is ensured.
The results are applied to regularized estimation of an ergodic diffusion observed at high frequency.
\end{abstract}

%%% Introduction
\section{Introduction}

We are concerned here with moment convergence in regularized estimation under multiple and mixed-rates asymptotics.
Let us begin with some basic background of the present study.
Suppose that we observe data whose distribution is indexed by a finite-dimensional parameter $\theta\in\Theta\subset\mathbb{R}^{\mathsf{p}}$.
Let $\theta_{0}\in\Theta$ denote the true value of $\theta$.
In order to estimate $\tz$ by means of the M-estimation theory, we introduce an appropriate (quasi-)likelihood or contrast function $\mbbh_{n}: \Omega\times\Theta\to\mathbb{R}$, and estimate an optimal parameter value by any point $\tes\in\argmin\mbbh_{n}$. For assessing asymptotic performance of $\hat{\theta}_{n}$ quantitatively, we look at the statistical random fields
\begin{align}
\mathbb{M}_{n}(u;\theta_{0})=\mathbb{H}_{n}(\theta_{0}+A_{n}(\theta_{0})u)-\mathbb{H}_{n}(\theta_{0}),\qquad u\in\mathbb{R}^{\mathsf{p}},
\label{hm:intro1}
\end{align}
%$\theta=(\al,\beta)\in\Theta\subset\mathbb{R}^{p+q}$, where the components $\al$ and $\beta$ will be assumed to be estimated at different rates of convergence. Let $\theta_{0}=(\al_{0},\beta_{0})\in\Theta$ denote the true value of $\theta$.
%In order to estimate $\tz$ by means of the M-estimation theory (e.g. \cite[Chapter 5]{van98}), we introduce an appropriate (quasi-)likelihood or contrast function $\mbbh_{n}: \Omega\times\Theta\to\mathbb{R}$, and estimate an optimal parameter value by any point $\hat{\theta}_{n}=(\hat{\al}_{n},\hat{\beta}_{n})\in\argmin\mbbh_{n}$. For assessing asymptotic performance of $\hat{\theta}_{n}$ quantitatively, we look at the statistical random fields
%\begin{align}
%\mathbb{M}_{n}(u;\theta_{0})=\mathbb{H}_{n}(\theta_{0}+A_{n}(\theta_{0})u)-\mathbb{H}_{n}(\theta_{0}),\qquad u\in\mathbb{R}^{p+q},
%\label{hm:intro1}
%\end{align}
where $A_{n}(\theta_{0})$ denotes the rate matrix such that $|A_{n}(\theta_{0})|\to 0$ as $n\to\infty$ and that the components may decrease at different rates, say, of multiple-rates asymptotic;
estimation under multiple-rates asymptotic has appeared in the literature of, for example, econometrics \cite{AntRen12}.
As is well-known, the weak convergence of $\mathbb{M}_{n}$ to some $\mathbb{M}_{0}$ over compact sets, the identifiability condition on $\mathbb{M}_{0}$, and the tightness of the scaled estimator $\hat{u}_{n}:=A_{n}(\theta_{0})^{-1}(\hat{\theta}_{n}-\theta_{0})$ make the ``argmin'' functional for $\mathbb{M}_{n}$ continuous: $\hat{u}_{n}\in\argmin\mathbb{M}_{n}\xrightarrow{\mathcal{L}}\argmin\mathbb{M}_{0}$.
In many statistical models, this procedure is enough to completely clarify the asymptotic distribution of all the components of $\hat{\theta}_{n}$.
We will refer this case to ``standard asymptotics".
See \cite[Chapter 5]{van98} for a general yet concise account of the M-estimation theory.

Beyond the weak convergence, convergence of moments of $\hat{u}_{n}$ serves as a fundamental tool when analyzing asymptotic behavior of the expectations of statistics depending on the estimator such as asymptotic bias and mean squared prediction error.
We refer to \cite{ChaIng11}, \cite{FinWei02}, \cite{IngYan14}, \cite{SakYos04}, \cite{UchYos01}, and \cite{UchYos06} for some related results. Also, \cite{AfeMar15-2} recently discussed optimal selection of random and $k$-fold cross-validation estimators, the theoretical backbone of which involves some moment bounds of the estimators used; the related paper \cite{AfeMar15} studied the uniform integrability of the ordinary least-squares estimator in the linear regression setting. In the standard asymptotics, the polynomial type large deviation inequality (PLDI) of \cite{Yos11}, which estimates the tail of $\mathcal{L}(\hat{u}_{n})$ in such a way that
\begin{align}
\sup_{r>0}\sup_{n>0}r^{L}P(\vert \hat{u}_{n}\vert\geq r)<\infty
\nn
%\label{eq:introp}
\end{align}
for a given $L>0$, provides us with a widely-applicable tool for verifying the convergence of moments.
In \cite{Yos11}, it has been shown that the PLDI can be proved under modest conditions especially when $\mathbb{M}_{n}$ admits a locally asymptotically quadratic (LAQ) structure\footnote{The sign in front of the quadratic term $(1/2)\Gamma_{0}(\theta_{0})[u,u]$ is different from the original LAQ of \cite{Yos11} since we consider minimization of (\ref{hm:intro1}).}
\begin{align}
\mathbb{M}_{n}(u;\theta_{0})=\Delta_{n}(\theta_{0})[u]+\frac{1}{2}\Gamma_{0}(\theta_{0})[u,u]+r_{n}(u;\theta_{0}), \qquad u\in\mathbb{R}^{\mathsf{p}},
\nonumber
\end{align}
where $\D_{n}(\tz)$ and $\Gam_{0}(\tz)$ are random linear and bilinear forms, respectively (we use the multi-linear-from notations: $A[u]:=\sum_{i}A_{i}u_{i}$ and $B[u,u]:=\sum_{i,j}B_{ij}u_{i}u_{j}$, and so on), and where $r_{n}(u;\tz)$ is expected to be stochastically ignorable for each $u$.
The LAQ property is known to be satisfied for many situations including non-ergodic (asymptotically mixed-normal type) models under multi-scaling.
A more detailed description is given in Appendix \ref{appendix-A}.

\medskip

In principle, any $M$-estimation procedure may have its ``regularized'' counterpart.
The logic of the sparse and more generally shrinkage estimation would be most clearly described by the context of multiple linear regression, with many deep theoretical interpretation such as geometrical (projection) characterization, variable selection, stabilized prediction performance, etc.; see e.g. \cite{BicLi06}, \cite[Chapter 3]{HasTibFri09}, and \cite{PotLee09}.
Statistical random fields associated with regularized estimation may not be LAQ, if the random field is of the ``mixed-rates" type in the sense of \cite{Rad08} where the target statistical random fields may have components converging at different rates.
In that case, it may even happen that the random function $\mbbm_{n}(u;\tz)$ diverges in probability for each $u$;
indeed, the sparse-type estimation, which is a particular case of regularized estimation, falls into the region of mixed-rates asymptotics.
In this case, convergence of moments does not follow from a direct application of \cite{Yos11}, and we are aware of only the following previous studies in this direction: 
\cite{UmeShiMasNin15} deduced the convergence of moments of a regularized sparse maximum-likelihood estimator of the generalized linear model, and applied it to verify the AIC type variable selection;
\cite{Shi16} (also \cite{MasShi14}) deduced the moment convergence in regularized estimation of a linear regression model with general regularization term.
Nevertheless, the proofs of these results made particular use of the special structure of the considered models and/or the convexity argument, and both of them do not tell us much about the case of general regularized $M$-estimation.

The primary contribution of this paper is to compose a first systematic study of not only asymptotic distribution but also convergence of moments associated with a class of regularized LAQ models under multiple and mixed-rates asymptotics.
In particular, in Theorem \ref{thm:pldi}, we will show how the PLDI of \cite{Yos11} can carry over to the case of mixed-rates $M$-estimation.

\medskip

This paper is organized as follows. Section \ref{hm:sec_setup} describes our model setup, based on which a series of basic asymptotic statements are given in Section \ref{hm:sec_asymptotics}.
In Section \ref{hm:sec_s.asymp} we will deal with the sparse asymptotics, and then the case of standard asymptotics is briefly mentioned in Section \ref{hm:sec_st.asymp.}. 
In Section \ref{hm:sec_tp.selection}, we will discuss a simple specific example of component-wise tuning-parameter choice. 
%The shrinkage effect is still useful for dependent-data models: it can effectively diminish non-significant factor involved in the model, resulting in a model-complexity assessment and/or selection.
In Section \ref{hm:sec_ergo.diff} we will apply the foregoing results to regularized estimation of an ergodic diffusion process observed at high frequency. The model is described by the Wiener-driven stochastic differential equation
\begin{align}
dX_{t}=a(X_{t},\alpha)dw_{t}+b(X_{t},\beta)dt,
\nonumber
\end{align}
and is known to be of multiple-scaling type. In the literature, \cite{DeGIac12} previously studied an adaptive-lasso type regularized estimation of the same model, and deduced the oracle property. 
In this study, we will derive the sparse consistency of the zero parameters, the asymptotic normality of the non-zero parameters, and the PLDI for a regularized estimator without resorting to convexity at all.

It is important to study associated practical computational algorithms, e.g. tuning parameter selection by making use of AIC (see Section \ref{hm:sec_aic}), and numerical optimization subjects such as the alternating direction method of multipliers. They are beyond the scope of this paper and we would like to leave it as one of future works.

%%%%%
%%%%%
\section{Setup}\label{hm:sec_setup}

Let us begin with description of the basic model setup for Section \ref{hm:sec_asymptotics}.
Throughout we are given an underlying probability space $(\Omega,\mcf,P)$. For the purpose of accelerating estimation performance, we consider $M$-estimation of an additive regularization type. We will focus on the case of two-scaling, where the target statistical parameter $\theta\in\Theta$, is divided into two parts:
\begin{align}
\theta=(\al,\beta).
\nonumber
\end{align}
We set $\al\in\mbbr^{p}$ and $\beta\in\mbbr^{q}$, and $\Theta=\Theta_{\alpha}\times\Theta_{\beta}$ to be a bounded convex domain in $\mbbr^{p+q}$.

We are given a function $M_{n}: \Omega\times\Theta\to\mbbr$, and regularization (possibly random) functions $\overline{R}_{n}^{a}(\al)$ and $\overline{R}_{n}^{b}(\beta)$. 
We then consider a contrast function $\mbbh_{n}: \Omega\times\Theta\to\mbbr$ of the form
\begin{align}
\mathbb{H}_{n}(\theta)=\mathbb{H}_{n}(\alpha,\beta)=M_{n}(\al,\beta)+\overline{R}_{n}^{a}(\al) +\overline{R}_{n}^{b}(\beta).
\label{hm:Hn}
\end{align}
The associated regularized $M$-estimator is defined to be any element (for brevity, implicitly assumed to exist)
\begin{align}
\tes \in \argmin_{\theta\in\overline{\Theta}}\mbbh_{n}(\theta).
\nonumber
\end{align}
We quantitatively distinguish zero parameters from non-zero ones. 
We denote by $\theta_{0}=(\al_{0},\beta_{0})$ the value we want to estimate (typically the true value of $\theta$) and assume that it takes the form $\al_{0}=(\al_0^{\circ},\al_0^{\ast})=((\al^{\circ}_{0,k'})_{k'},(\al^{\ast}_{0,k''})_{k''})$ and $\beta_{0}=(\beta_0^{\circ},\beta_0^{\ast})=((\beta^{\circ}_{0,l'})_{l'},(\beta^{\ast}_{0,l''})_{l''})$ with
\begin{align}
\al^{\circ}_{0,k'}=0,\quad \beta^{\circ}_{0,l'}=0, \quad \al^{\ast}_{0,k''}\ne 0, \quad \beta^{\ast}_{0,l''}\ne 0.
\nonumber
\end{align}
We set $\al^{\circ}_{0}\in\mbbr^{p^{\circ}}$, $\beta^{\circ}_{0}\in\mbbr^{q^{\circ}}$, $\al^{\ast}_{0}\in\mbbr^{p^{\ast}}$ and $\beta^{\ast}_{0}\in\mbbr^{q^{\ast}}$ with $p^{\circ},q^{\circ},p^{\ast},q^{\ast}\in\mbbn$; then, $p=p^{\circ}+p^{\ast}$ and $q=q^{\circ}+q^{\ast}$. 
Correspondingly, we write $\theta=(\theta^{\circ},\theta^{\ast})$ with $\theta^{\circ}=(\alpha^{\circ},\beta^{\circ})$ and $\theta^{\ast}=(\alpha^{\ast},\beta^{\ast})$. 
We also write $\hat{\theta}_{n}=(\hat{\alpha}_{n},\hat{\beta}_{n})=(\hat{\alpha}_{n}^{\circ},\hat{\alpha}_{n}^{\ast},\hat{\beta}_{n}^{\circ},\hat{\beta}_{n}^{\ast})$ 
with $\tes^{\circ}=(\aes^{\circ},\bes^{\circ})$ and $\tes^{\ast}=(\aes^{\ast},\bes^{\ast})$. 
For clarity we focus on regularization terms of the form
\begin{align}
\overline{R}_{n}^{a}(\al) = \sum_{k=1}^{p}\lam^{a}_{n,k}R^{a}(\al_{k}), \qquad \overline{R}_{n}^{b}(\beta) = \sum_{l=1}^{q}\lam^{b}_{n,l}R^{b}(\beta_{l}).
\label{hm:regu.term.add1}
\end{align}
The regularization-function form subsumes many existing ones, although not essential for our basic asymptotic results given in Sections  \ref{hm:sec_s.asymp} and \ref{hm:sec_st.asymp.}. For later reference we write more specifically
\begin{align}
\overline{R}_{n}^{a}(\al) 
&=\overline{R}_{n}^{a \circ}(\al^{\circ}) + \overline{R}_{n}^{a \ast}(\al^{\ast}) 
= \sum_{k'=1}^{p^{\circ}}\lam^{a \circ}_{n,k'}R^{a}(\al^{\circ}_{k'}) + \sum_{k''=1}^{p^{\ast}}\lam^{a \ast}_{n,k''}R^{a}(\al^{\ast}_{k''}),
\label{hm:Ra_form}\\
\overline{R}_{n}^{b}(\beta) 
&=\overline{R}_{n}^{b \circ}(\beta^{\circ}) + \overline{R}_{n}^{b \ast}(\beta^{\ast}) 
= \sum_{l'=1}^{q^{\circ}}\lam^{b \circ}_{n,l'}R^{b}(\beta^{\circ}_{l'}) + \sum_{l''=1}^{q^{\ast}}\lam^{b \ast}_{n,l''}R^{b}(\beta^{\ast}_{l''}),
\label{hm:Rb_form}
\end{align}
where:
\begin{itemize}
\item $\lam^{a \circ}_{n,k'}$, $\lam^{a \ast}_{n,k''}$,  $\lam^{b \circ}_{n,l'}$ and $\lam^{b \ast}_{n,l''}$ are non-negative random variables;
\item $R^{a}(\cdot)$, $R^{b}(\cdot)$ are nonrandom non-negative functions on $\mbbr$ such that $R^{a}(0)=R^{b}(0)=0$;
\item For all $a_{0},\ b_{0}\not=0$ and $k>0$, there exists a constant $C=C(a_{0},b_{0},k)>0$ such that
\begin{align}
\sup_{(a',b'):\,|a'|\vee|b'|\leq k}\frac{|R^{a}(a')-R^{a}(a_{0})|+|R^{b}(b')-R^{b}(b_{0})|}{|a'-a_{0}|+|b'-b_{0}|}\leq C. 
\label{hm:llc}
\end{align}
\end{itemize}
The last condition (local Lipschitz continuity) is a technical one. Further conditions on $M_{n}$, $\overline{R}_{n}^{a}(\al)$ and $\overline{R}_{n}^{b}(\beta)$ will be imposed later on; in Section \ref{hm:sec_tp.selection}, we will briefly discuss about how to set the regularization terms in some specific ways.

We will deal with a situation where the non-zero part of the first component $\al$ can be estimated faster than that of the second component $\beta$; 
more specifically, we will suppose that the sequence
\begin{align}
\left( s_{n}^{-1}(\aes^{\ast}-\al_{0}^{\ast}),\, t_{n}^{-1}(\bes^{\ast}-\beta_{0}^{\ast}) \right)
\nonumber
\end{align}
has a non-trivial asymptotic distribution for some possibly different positive sequence $(s_{n})$ and $(t_{n})$, both tending to zero and satisfying that $s_{n}=o(t_{n})$. Although not explicitly mentioned, we presuppose that the ``principal'' part $M_{n}(\theta)$ reasonably makes sense even without regularization terms 
$\overline{R}_{n}^{a}(\al)+\overline{R}_{n}^{b}(\beta)$; most typically, the un-regularized case, where $\mathbb{H}_{n}(\theta)=M_{n}(\theta)$, corresponds to a negative of a (quasi) log-likelihood.
%A concrete example is the regularized Gaussian quasi-likelihood estimation of ergodic diffusion observed at high frequency, which we will discuss in Section \ref{hm:sec_ergo.diff}.
We should note that the additive regularization can be interpreted as incorporating a prior information about the parameter of interest; see Section \ref{hm:sec_tp.selection}.

We end this section with noting that an extension to cases of more-than-two scaling is a trivial matter while making notation messy,
and also that we can consider the usual single-scaling case as well simply by omitting the parameter $\beta$.

%%%%%
%%%%%
\section{Asymptotics}\label{hm:sec_asymptotics}

\subsection{Sparse case}\label{hm:sec_s.asymp}
Under the setting described in Section \ref{hm:sec_setup}, we consider regularity conditions under which the following properties hold without assuming the convexity of $\mbbh_{n}$.
\begin{enumerate}
\item The (weak) consistency of $\tes=(\tes^{\circ},\tes^{\ast})$.
\item The asymptotic distributions:
\begin{enumerate}
\item The sparse consistency of $\tes^{\circ}$, i.e. $P(\tes^{\circ}=0)\to 1$;
\item The asymptotic distribution of $\tes^{\ast}$ at possibly multiple rates of convergence (via a matrix norming).
\end{enumerate}
\item The uniform tail-probability estimate of $\tes=(\tes^{\circ},\tes^{\ast})$.
\end{enumerate}
Hereafter, we will prove the results in the above order: Theorem \ref{hm:thm_consistency} describes the consistency of $\hat{\theta}_{n}$. In Theorem \ref{hm:thm_rate} we derive the auxiliary convergence rate, followed by the sparse consistency of $\tes^{\circ}$ in Theorem \ref{ys:thm_sc}. We also derive the asymptotic distribution of $\hat{\theta}_{n}^{\ast}$ in Theorem \ref{hm:thm_AL}. Finally, we will present the uniform tail-probability estimate of $\tes=(\tes^{\circ},\tes^{\ast})$ in Theorem \ref{thm:pldi}.

%%%
\subsubsection{Consistency}\label{hm:sec_consistency}

We impose the uniform law of large numbers plus identifiability condition, and additionally some stochastic-order conditions on the regularization terms. 
Recall that the parameter space $\Theta$ is a bounded convex domain.

\begin{assum}\label{hm:A_c}\ 
\begin{enumerate}
\item $(s_{n})$ and $(t_{n})$ are positive nonrandom sequences such that $\max(s_{n}, t_{n})\to 0$ and that $s_{n}=o(t_{n})$.
\item There exist continuous random functions $\overline{M}^{a}_{0}: \Omega\times\Theta_{\al} \to \mbbr$ and $\overline{M}^{b}_{0}: \Omega\times\Theta \to \mbbr$ such that:
\begin{enumerate}
\item $\displaystyle \sup_{\al}\left|s_{n}^{2}\left\{M_{n}(\al,\beta_{0}) - M_{n}(\al_{0},\beta_{0})\right\} 
- \overline{M}_{0}^{a}(\al)\right|$ \\
$\displaystyle +\sup_{\theta}\left|t_{n}^{2}\left\{M_{n}(\al,\beta) - M_{n}(\al,\beta_{0})\right\} 
- \overline{M}_{0}^{b}(\theta)\right| \cip 0$;
\item $\ds{\argmin_{\al}\overline{M}_{0}^{a}(\al)=\{\al_{0}\}}$ a.s. 
and $\ds{\argmin_{\beta}\overline{M}_{0}^{b}(\al_{0},\beta)=\{\beta_{0}\}}$ a.s.
\end{enumerate}
\item $\ds{\sup_{\al}\left|s_{n}^{2}\overline{R}_{n}^{a}(\al)\right| + \sup_{\beta}\left| t_{n}^{2}\overline{R}_{n}^{b}(\beta)\right| \cip 0}$.
\label{ys:consistency-3}
\end{enumerate}
\end{assum}

We will take advantage of a general result concerning mixed-rates asymptotics.

\begin{lem}\label{Rad08_jal.thm}
Let $\mbbm_{n}(u,v): \mbbr^{p}\times\mbbr^{q}\to\mbbr$ be a random function, and $(\hat{u}_{n},\hat{v}_{n})$ and $(\hat{u}_{0},\hat{v}_{0})$ random variables taking their values in $\mbbr^{p}\times\mbbr^{q}$. Assume that $\mbbm_{n}$ is of the form
\begin{align}
\mbbm_{n}(u,v) = \bar{a}_{n}f_{n}(u) + \bar{b}_{n}g_{n}(u,v),\qquad (u,v)\in\mathbb{R}^{p}\times\mathbb{R}^{q},
\nn
%\label{hm:Rad.Thm.M}
\end{align}
where:
\begin{itemize}
\item[(L1)] $\bar{a}_{n}$ and $\bar{b}_{n}$ are positive constants such that $\bar{b}_{n}=o(\bar{a}_{n})$;
\item[(L2)] $(f_{n}(\cdot),g_{n}(\cdot,\cdot))\stackrel{\mathcal{L}}{\to}(f_{0}(\cdot),g_{0}(\cdot,\cdot))$ in $\mathcal{C}(K_{f}\times K_{g})$ for every compact 
$K_{f}\times K_{g}\subset\mathbb{R}^{p} \times \mathbb{R}^{q}$;
\item[(L3)] $\mbbm_{n}(\hat{u}_{n},\hat{v}_{n}) \le \inf_{(u,v)}\mbbm_{n}(u,v) + o_{p}(\bar{b}_{n})$;
\item[(L4)] $(\hat{u}_{n},\hat{v}_{n}) = O_{p}(1)$;
\item[(L5)] $u\mapsto f_{0}(u)$ has a.s. unique minimum at $u=\hat{u}_{0}$;
\item[(L6)] $v\mapsto g_{0}(\hat{u}_{0},v)$ has a.s. unique minimum at $v=\hat{v}_{0}$.
\end{itemize}
Then we have the following.
\begin{enumerate}
\item $(\hat{u}_{n},\hat{v}_{n}) \stackrel{\mathcal{L}}{\to}(\hat{u}_{0},\hat{v}_{0})$ under the conditions (L1) to (L6).
\item $\hat{u}_{n} \stackrel{\mathcal{L}}{\to} \hat{u}_{0}$ under the conditions (L1) to (L5).
\end{enumerate}
\end{lem}

\begin{proof}
The first claim is just a simplified version of \cite[Theorem 1]{Rad08}. The second one follows on applying the usual argmax theorem to the rescaled function $u\mapsto \bar{a}_{n}^{-1}\mbbm_{n}(u,\hat{v}_{n})=f_{n}(u)+(\bar{b}_{n}/\bar{a}_{n})g_{n}(u,\hat{v}_{n})$:
this random function admits an approximate minimizer $\hat{u}_{n}$ and weakly converges to $f_{0}$ in $\mathcal{C}(K_{f})$ for every compact $K_{f}\subset\mbbr^{p}$.
\end{proof}

\begin{rem}{\rm 
Although we do not use the second claim of Lemma \ref{Rad08_jal.thm} in this paper, it may be useful when considering stepwise estimation.
%where there is an original contrast (or quasi-likelihood) function but some step-by-step strategy is taken for estimating parameter components that can be estimated more quickly than the others.
See \cite{UchYos12} for an ergodic diffusion model.
}\end{rem}

To apply the first claim of Lemma \ref{Rad08_jal.thm} under Assumption \ref{hm:A_c} with $(u,v)=(\al,\beta)$, we set $\mbbm_{n}(\al,\beta)=\mbbh_{n}(\al,\beta)-\mbbh_{n}(\al_{0},\beta_{0})$, $\bar{a}_{n}=s_{n}^{-2}$, $\bar{b}_{n}=t_{n}^{-2}$:
\begin{align}
f_{n}(\al) &= s_{n}^{2}\left( M_{n}(\al,\beta_{0}) - M_{n}(\al_{0},\beta_{0}) + \overline{R}_{n}^{a}(\al)-\overline{R}_{n}^{a}(\al_{0})\right), \nn\\
g_{n}(\al,\beta) &= t_{n}^{2}\left( M_{n}(\al,\beta) - M_{n}(\al,\beta_{0}) + \overline{R}_{n}^{b}(\beta)-\overline{R}_{n}^{b}(\beta_{0})\right). \nn
\nonumber
\end{align}
According to the a.s. continuity of the random functions $\al\mapsto\overline{M}^{a}_{0}(\al)$ and $\theta\mapsto\overline{M}^{b}_{0}(\theta)$, it is straightforward to verify all the conditions in Lemma \ref{Rad08_jal.thm} under Assumption \ref{hm:A_c}, hence the following claim:

\begin{thm}\label{hm:thm_consistency}
We have $\tes\cip\tz$ under Assumption \ref{hm:A_c}.
\end{thm}

%%%
\subsubsection{Rates of convergence}\label{hm:sec_rc}

Next we prove that $(\hat{u}_{n}, \hat{v}_{n})=O_{p}(1)$ where
\begin{align}
\hat{u}_{n} := s_{n}^{-1}(\aes-\al_{0}), \qquad \hat{v}_{n} := t_{n}^{-1}(\bes-\beta_{0})
\label{hm:uhat.vhat.def}
\end{align}
under additional conditions. To this end we introduce the following general result, which is a simplified version of \cite[Lemma 1]{Rad08}, to deduce ``correct'' convergence rate of $\tes^{\ast}$ and a ``auxiliary'' convergence rate of $\tes^{\circ}$. Here, the auxiliary convergence rate may not be the optimal convergence rate for $\tes^{\circ}$, for which we will show more precisely the sparse consistency. Nevertheless, the auxiliary rate should be useful for verification of the conditions for the sparse consistency of $\tes^{\circ}$ (see Remark \ref{hm:rem_bridge} below). Let $[a]_{+}:=\max(a,0)$ for $a\in\mbbr$.

\begin{lem}\label{Rad08_lem1-}
Let $\xi$ denote either $\al$ or $\beta$, and 
assume that the real-valued random function $\overline{\mbbh}_{n}(\xi)$ satisfies the following conditions:
\begin{enumerate}
\item $\overline{\mbbh}_{n}(\xes) \le \overline{\mbbh}_{n}(\xi_{0})$ a.s.
%\item $\xes\cip\xi_{0}$.
\item There exist random functions $U_{n}(\xi)$ and $V_{n}(\xi)$ such that
\begin{align}
\overline{\mbbh}_{n}(\xi) - \overline{\mbbh}_{n}(\xi_{0}) = U_{n}(\xi) - V_{n}(\xi)
\nonumber
\end{align}
where, for some random variable $U_{0}>0$ a.s., constants $0\le \gam<\rho$, and positive nonrandom sequence $(k_{n})$ such that $k_{n}\to 0$, we have:
\begin{enumerate}
\item $\ds{P\left( U_{n}(\xes) \ge |\xes-\xi_{0}|^{\rho}U_{0}\right) \to 1}$;
\item $\ds{[V_{n}(\xes)]_{+}=O_{p}(k_{n}|\xes-\xi_{0}|^{\gam})}$.
\end{enumerate}
\end{enumerate}
Then $k_{n}^{-1/(\rho-\gam)}(\xes-\xi_{0})=O_{p}(1)$.
\end{lem}

To deduce Lemma \ref{Rad08_lem1-}, observe that on the set $\{|\xes-\xi_{0}|^{\rho}\le U_{0}^{-1}U_{n}(\xes)\}$ we have $|\xes-\xi_{0}|^{\rho} \le U_{0}^{-1}\{V_{n}(\xes)+\overline{\mbbh}_{n}(\xes) - \overline{\mbbh}_{n}(\xi_{0})\} \le U_{0}^{-1}V_{n}(\xes) \le U_{0}^{-1}[V_{n}(\xes)]_{+}=L_{n}k_{n}|\xes-\xi_{0}|^{\gam}$, where $L_{n}=O_{p}(1)$. By applying \cite[Lemma 3]{Rad08}, for any $\del>0$ we can find a tight sequence $(N_{n})$ in $\mbbr$ for which $L_{n}k_{n}|\xes-\xi_{0}|^{\gam} \le \del |\xes-\xi_{0}|^{\rho}+N_{n}k_{n}^{\rho/(\rho-\gam)}$. Picking a $\del<1$, we conclude that $|\xes-\xi_{0}|\le k_{n}^{1/(\rho-\gam)}N_{n}$, hence the assertion.

\medskip

In what follows, for any matrix $A$ we write $A^{\otimes2}=AA^{\top}$ with $\top$ denoting the transpose,
and for any real symmetric matrix $A$ we denote by $\lam_{\min}(A)$ the smallest eigenvalue of $A$. Coming back to our model, we next impose:

\begin{assum}\label{hm:A_rc}\ 
\begin{enumerate}
\item $M_{n} \in \mcc^{3}(\Theta)$ a.s., and it holds that:
\label{ys:R_cond.1}
\begin{enumerate}
\item $\ds{\sup_{\beta}\left|s_{n}\p_{\al}M_{n}(\al_{0},\beta)\right| + \left|t_{n}\p_{\beta}M_{n}(\tz)\right| = O_{p}(1)}$;
\label{ys:R_cond.1-a}
\item $\ds{\sup_{\al}\left|s_{n}t_{n}\p_{\al}\p_{\beta}M_{n}(\al,\beta_{0})\right| = O_{p}(1)}$;
\label{ys:R_cond.1-b}
\item $\ds{\sup_{\theta}\left|s_{n}^{2}\p_{\zeta}\p_{\alpha}^{2}M_{n}(\theta)\right| + \sup_{\theta}\left|t_{n}^{2}\p_{\zeta}\p_{\beta}^{2}M_{n}(\theta)\right| = O_{p}(1)}$ for $\zeta=\al, \beta$;
\label{ys:R_cond.1-c}
\item \label{ys:R_cond.1-d}There exist symmetric random functions $\Gam_{0}^{\al}: \Omega\times\Theta_{\al} \to \mbbr^{p}\otimes\mbbr^{p}$ and 
$\Gam_{0}^{\beta}: \Omega\times\Theta \to \mbbr^{q}\otimes\mbbr^{q}$ such that
\begin{align}
\left|s_{n}^{2}\p_{\al}^{2}M_{n}(\tz) - \Gam_{0}^{\al}(\al_{0})\right|+\left|t_{n}^{2}\p_{\beta}^{2}M_{n}(\tz) - \Gam_{0}^{\beta}(\tz)\right|\cip 0,
\nonumber
\end{align}
with $\ds{\lam_{\min}\big( \Gam_{0}^{\al}(\al_{0}) \big)\wedge\lam_{\min}\big( \Gam_{0}^{\beta}(\theta_{0})\big) > 0}$ a.s.
\end{enumerate}
\item $s_{n}\lam^{a \ast}_{n,k''}=O_{p}(1)$ and $t_{n}\lam^{b \ast}_{n,l''}=O_{p}(1)$ for each $k''$ and $l''$.
\label{ys:R_cond.3}
\end{enumerate}
\end{assum}

\begin{rem}
\upshape
Note that Assumption \ref{hm:A_rc}.\ref{ys:R_cond.3} is only concerned with the non-zero parameter parts, which of course are unknown \textit{a priori} in practice. Hence, as is well recognized as a common situation in adaptive type sparse estimation, some appropriate data-driven choices of the weights are desirable. There would be many possibilities for this issue. See Section \ref{hm:sec_tp.selection} for more discussions.
\end{rem}

\medskip

Let Assumptions \ref{hm:A_c} and \ref{hm:A_rc} hold, so that $\tes\cip\tz$ by Theorem \ref{hm:thm_consistency}. We apply Lemma \ref{Rad08_lem1-} for proving $\hat{u}_{n}=O_{p}(1)$ and $\hat{v}_{n}=O_{p}(1)$ separately. To show $\hat{u}_{n}=O_{p}(1)$, set
\begin{align}
\overline{\mbbh}_{n}(\al) = s_{n}^{2}\mbbh_{n}(\al,\bes).
\nonumber
\end{align}
For notational convenience, for any random function $F_{n}(\theta)$ we will write $F_{n}(\theta)=\ulop{1}$ and $F_{n}(\theta)=\usop{1}$ if $\sup_{\theta}|F_{n}(\theta)|=O_{p}(1)$ and $\sup_{\theta}|F_{n}(\theta)|=o_{p}(1)$, respectively. The first condition in Lemma \ref{Rad08_lem1-} is trivial. To deduce the second one, making use of the third-order Taylor expansion we derive $\overline{\mbbh}_{n}(\aes)-\overline{\mbbh}_{n}(\al_{0})=U_{n}(\aes) - V_{n}(\aes)$ for
\begin{align}
U_{n}(\aes) &:=\frac{1}{2}s_{n}^{2}\p_{\al}^{2}M_{n}(\tilde{\al}_{n},\bes) \left[(\aes-\al_{0})^{\otimes 2}\right]
+ s_{n}^{2}\sum_{k'=1}^{p^{\circ}}\lam^{a \circ}_{n,k'}R^{a}(\hat{\al}^{\circ}_{n,k'})
\nn \\
&=\frac{1}{2}\left(s_{n}^{2}\p_{\al}^{2}M_{n}(\al_{0},\bes) + s_{n}^{2}\p_{\al}^{3}M_{n}(\check{\al}_{n},\bes)[\tilde{\al}_{n}-\al_{0}]\right)\left[(\hat{\al}_{n}-\al_{0})^{\otimes 2}\right] + s_{n}^{2}\sum_{k'=1}^{p^{\circ}}\lam^{a \circ}_{n,k'}R^{a}(\hat{\al}^{\circ}_{n,k'}) \nn\\
&=\frac{1}{2}\left\{\Gam^{\al}_{0}(\al_{0}) + \ulop{|\aes-\al_{0}|\vee|\bes-\beta_{0}|}\right\}\left[(\hat{\al}_{n}-\al_{0})^{\otimes 2}\right] + s_{n}^{2}\sum_{k'=1}^{p^{\circ}}\lam^{a \circ}_{n,k'}R^{a}(\hat{\al}^{\circ}_{n,k'}) \nn\\
&=\frac{1}{2}\left\{\Gam^{\al}_{0}(\al_{0}) + \usop{1}\right\}\left[(\hat{\al}_{n}-\al_{0})^{\otimes 2}\right] + s_{n}^{2}\sum_{k'=1}^{p^{\circ}}\lam^{a \circ}_{n,k'}R^{a}(\hat{\al}^{\circ}_{n,k'}),
\label{hm:Ua}\\
V_{n}(\hat{\al}_{n}) &:= -s_{n}^{2}\p_{\al}M_{n}(\al_{0}, \bes) [\hat{\al}_{n}-\al_{0}] 
- s_{n}^{2}\sum_{k''=1}^{p^{\ast}}\lam^{a \ast}_{n,k''}\left( R^{a}(\hat{\al}^{\ast}_{n,k''})-R^{a}(\al^{\ast}_{0,k''})\right),
\nn
%\label{hm:Va}
\end{align}
where the points $\tilde{\al}_{n}$ and $\check{\al}_{n}$ are located on the segments connecting $\al_{0}$ and $\aes$, and $\al_{0}$ and $\tilde{\al}_{n}$, respectively. Note that the non-negativity of $R^{a}$ enables us to ignore the second term of the right-hand side of \eqref{hm:Ua} when estimating $U_{n}(\hat{\al}_{n})$ from below. 
Also, under the local Lipschitz continuity \eqref{hm:llc} of $R^{a}$, we see that the conditions of Lemma \ref{Rad08_lem1-} are satisfied with $\gam=1$, $\rho=2$, (e.g.) $U_{0}=\lam_{\min}\left( \Gam_{0}^{\al}(\al_{0}) \right)/4$, and $k_{n}=s_{n}$: we have $U_{n}(\aes) \ge (1/2)\{o_{p}(1)+\lam_{\min}(\Gam_{0}^{\al}(\al_{0}))\}|\aes-\al_{0}|^{2}$ and $|V_{n}(\aes)|\le O_{p}(s_{n}|\aes-\al_{0}|)$. Hence $\hat{u}_{n}=O_{p}(1)$ is proved. 
To deduce $\hat{v}_{n}=O_{p}(1)$, we can follow exactly the same way as above along with
\begin{align}
\overline{\mbbh}_{n}(\beta) = t_{n}^{2}\left\{\mbbh_{n}(\aes,\beta)-\mbbh_{n}(\aes,\beta_{0})\right\}
\nonumber
\end{align}
in place of $\overline{\mbbh}_{n}(\al) = s_{n}^{2}\mbbh_{n}(\al,\bes)$.

\begin{thm}\label{hm:thm_rate}
We have $(\hat{u}_{n}, \hat{v}_{n})=O_{p}(1)$ under Assumptions \ref{hm:A_c} and \ref{hm:A_rc}.
\end{thm}

%\begin{rem}{\rm 
%Another way to establish the tightness of $\{(\hat{u}_{n},\hat{v}_{n})\}_{n}$ is to ``follow'' the derivation of the PLDI \cite[Theorem 1]{Yos11} without reference to the polynomial orders of the tail-probability estimate.
%}\end{rem}

%%%
\subsubsection{Sparse consistency}\label{hm:sec_sc}

The sparse consistency $P(\tes^{\circ}=0)\to 1$ implies that the asymptotic distribution of $\tes^{\circ}$ degenerates at the origin with arbitrarily fast rate of convergence: $r_{n}\hat{\theta}_{n}^{\circ}=o_{p}(1)$ for arbitrary $r_{n}\to\infty$. The next general result is a variant of \cite[Theorem 2]{Rad08}, which is a tailor-made tool to establish the property.

\begin{lem}\label{Rad08_thm2-}
Let $\xi$ denote either $\al$ or $\beta$ (so $\xi=(\xi^{\circ},\xi^{\ast})$ and $\xi_0=(0,\xi_0^{\ast})$), 
and assume that the real-valued random function $\overline{\mbbh}_{n}(\xi)=\overline{\mbbh}_{n}(\xi^{\circ},\xi^{\ast})$ satisfies the following conditions.
\begin{enumerate}
\item $\overline{\mbbh}_{n}(\xes^{\circ},\xes^{\ast}) \le \overline{\mbbh}_{n}(0,\xes^{\ast})$ a.s.
%\item $\xes^{\circ}\cip 0$.
\item There exist random functions $U_{n}(\xi)$ and $V_{n}(\xi)$ such that
\begin{align}
\overline{\mbbh}_{n}(\xi) - \overline{\mbbh}_{n}(0,\xi^{\ast}) = U_{n}(\xi) - V_{n}(\xi),
\nonumber
\end{align}
where it holds that for some random variable $U_{0}>0$ a.s. and constants $\rho>0$,
\begin{align}
P\bigg[ \left\{[V_{n}(\xes)]_{+}=0 \right\} \cap \left\{U_{n}(\xes) \ge |\hat{\xi}_{n}^{\circ}|^{\rho}U_{0}\right\}\bigg] \to 1.
\nn%\label{Rad08_thm2-.eq1}
\end{align}
\label{Rad08_thm2-.eq1}
\end{enumerate}
Then $P(\xes^{\circ}=0)\to 1$.
\end{lem}

Lemma \ref{Rad08_thm2-} follows on observing the following: 
on the event $A_{n}:=\{[V_{n}(\xes)]_{+}=0\}\cap\{ U_{n}(\xes^{\circ},\xes^{\ast}) \ge |\hat{\xi}_{n}^{\circ}|^{\rho}U_{0} \}$, we have $|\hat{\xi}_{n}^{\circ}|^{\rho} \le U_{0}^{-1}\{ V_{n}(\xes) + \overline{\mbbh}_{n}(\xes) - \overline{\mbbh}_{n}(0,\xes^{\ast}) \} \le U_{0}^{-1}[V_{n}(\xes)]_{+}=0$. Hence $P(|\hat{\xi}^{\circ}_{n}|=0)$ can be bounded below by $P(A_{n})\to 1$.

\begin{rem}{\rm 
To conclude $P(\xes^{\circ}=0)\to 1$ we may replace the second condition of Lemma \ref{Rad08_thm2-} by
\begin{align}
P\bigg[ \left\{[V_{n}(\xes)]_{+}=0 \right\} \cap \left\{U_{n}(\xes) \ge |\hat{\xi}_{n,m'}^{\circ}|^{\rho}U_{0}\right\}\bigg] \to 1
\nn
\end{align}
for each $m'$ running through $\{1,\dots,p^{\circ}\}$ or $\{1,\dots,q^{\circ}\}$ according as $\xi=\al$ or $\beta$; for each $m'$, the same proof as above leads to $P(\hat{\xi}_{n,m'}^{\circ}=0)\to 1$.
}\end{rem}

\medskip

Let us go back to our main context. We keep imposing Assumptions \ref{hm:A_c} and \ref{hm:A_rc}.
Denoting by $\Gam_{0}^{\al^{\circ}}(\al)$ (resp. $\Gam_{0}^{\beta^{\circ}}(\theta)$) the upper left $p^{\circ}\times p^{\circ}$ part of $\Gam_{0}^{\al}(\al)$ 
(resp. the upper left $q^{\circ}\times q^{\circ}$ part of $\Gam_{0}^{\beta}(\theta)$), we have
\begin{align}
\left|s_{n}^{2}\p_{\al^{\circ}}^{2}M_{n}(\tz) - \Gam_{0}^{\al^{\circ}}(\al_{0})\right|+\left|t_{n}^{2}\p_{\beta^{\circ}}^{2}M_{n}(\tz) - \Gam_{0}^{\beta^{\circ}}(\tz)\right|\cip 0
\label{hm:for.sparse.LLN}
\end{align}
with $\lam_{\min}\big(\Gam_{0}^{\al^{\circ}}(\al_{0}) \big) \wedge \lam_{\min}\big( \Gam_{0}^{\beta^{\circ}}(\tz)\big) > 0$ a.s. 
The next assumption imposes component-wise ``good balances'' between the quadratic expansions of $M_{n}$ and the regularization terms.

\begin{assum}\label{ys:sc}
There exist constants $\underline{a}_{k'}, \underline{b}_{l'}\in(0,1/2)$ such that
\begin{align}
& P\bigg(s_{n}^{2}\partial_{\alpha_{k'}^{\circ}}M_{n}(0,\dots,0,\hat{\al}^{\circ}_{n,k'},\dots,\hat{\al}^{\circ}_{n,p^{\circ}},
\aes^{\ast},\bes) \hat{\al}^{\circ}_{n,k'} \nn\\
&{}\qquad 
+s_{n}^{2}\lam^{a \circ}_{n,k'}R^{a}(\hat{\al}^{\circ}_{n,k'})\geq -\underline{a}_{k'}\lam_{\min}\big(\Gam_{0}^{\al^{\circ}}(\al_{0}) \big) |\hat{\al}^{\circ}_{n,k'}|^{2}\bigg)\rightarrow 1, \nn\\
& P\bigg(t_{n}^{2}\partial_{\beta_{l'}^{\circ}}M_{n}(\aes,0,\dots,0,\hat{\beta}^{\circ}_{n,l'},\dots,\hat{\beta}^{\circ}_{n,q^{\circ}},\hat{\beta}^{\ast}_{n}) \hat{\beta}_{n,l'}^{\circ} \nn\\
&{}\qquad + t_{n}^{2}\lam^{b \circ}_{n,l'}R^{b}(\hat{\beta}^{\circ}_{n,l'})\geq -\underline{b}_{l'}\lam_{\min}\big(\Gam_{0}^{\beta^{\circ}}(\tz) \big)
|\hat{\beta}_{n,l'}^{\circ}|^{2}\bigg)\rightarrow 1
\nonumber
\end{align}
for each $k'\in\{1,\dots,p^{\circ}\}$ and $l'\in\{1,\dots,q^{\circ}\}$.
%\begin{align}
%& P\left(s_{n}^{2}\partial_{\alpha^{\circ}}M_{n}(0,\aes^{\ast},\bes)[\aes^{\circ}]+s_{n}^{2}\overline{R}_{n}^{a \circ}(\aes^{\circ})\geq -a'\lam_{\min}\big(\Gam_{0}^{\alpha^{\circ}}(\al_{0})\big)|\aes^{\circ}|^{2}\right)\rightarrow 1, \nn\\
%& P\left(t_{n}^{2}\partial_{\beta^{\circ}}M_{n}(\aes,0,\bes^{\ast})[\bes^{\circ}]+t_{n}^{2}\overline{R}_{n}^{b \circ}(\bes^{\circ})\geq -b'\lam_{\min}\big(\Gam_{0}^{\beta^{\circ}}(\tz)\big)|\bes^{\circ}|^{2}\right)\rightarrow 1.
%\nonumber
%\end{align}
%%There exist random variables $H_{n}^{\al}$, $H_{n}^{\beta}$, $\epsilon^{\alpha}_{0}>0$, $\epsilon^{\beta}_{0}>0$ a.s. and constants $a', b'\in(0,1)$ such that:
%%\begin{enumerate}
%%\item $P(s_{n}^{2}\partial_{\alpha^{\circ}}M_{n}(0,\aes^{\ast},\bes)[\aes^{\circ}]+s_{n}^{2}\overline{R}_{n}^{a \circ}(\aes^{\circ})\geq H_{n}^{\alpha}|\aes^{\circ}|^{2})\rightarrow1$ and $P(t_{n}^{2}\partial_{\beta^{\circ}}M_{n}(\aes,0,\bes^{\ast})[\bes^{\circ}]+t_{n}^{2}\overline{R}_{n}^{b \circ}(\bes^{\circ})\geq H_{n}^{\beta}|\bes^{\circ}|^{2})\rightarrow1$; 
%%\item $P(H_{n}^{\al}+a'\inf_{\alpha}\lambda_{\min}(\Gamma^{\al^{\circ}}_{0}(\alpha))/2\geq\epsilon_{0}^{\alpha})\rightarrow1$ and $P(H_{n}^{\beta}+b'\inf_{\beta}\lambda_{\min}(\Gamma^{\beta^{\circ}}_{0}(\al_{0},\beta))/2\geq\epsilon_{0}^{\beta})\rightarrow1$. 
%%\end{enumerate}
\end{assum}

%\begin{rem}
%\upshape
%It is allowed that $H_{n}^{\alpha}\cip\infty$ and/or $H_{n}^{\beta}\cip\infty$, i.e. $P(H_{n}^{\alpha}>k)\rightarrow1$ and $P(H_{n}^{\beta}>k)\rightarrow1$ for any $k>0$. Such an example is given in \cite[Section 5]{Rad08}.
%\end{rem}

First, we apply Lemma \ref{Rad08_thm2-} with $\overline{\mbbh}_{n}(\alpha)=s_{n}^{2}\mbbh_{n}(\alpha,\bes)$ to deduce $P(\aes^{\circ}=0)\rightarrow1$. To show the second condition, by representing the $M_{n}$-part as a sum of $\al$-component-wise differences and then applying the second-term Taylor expansion, we have
\begin{align}
& \overline{\mbbh}_{n}(\aes)-\overline{\mbbh}_{n}(0,\aes^{\ast}) \nn\\
&=s_{n}^{2}\{\mbbh_{n}(\aes,\bes)-\mbbh_{n}(0,\aes^{\ast},\bes)\} \nn \\
&=s_{n}^{2}\{M_{n}(\aes,\bes)-M_{n}(0,\aes^{\ast},\bes)+\overline{R}_{n}^{a \circ}(\hat{\al}_{n}^{\circ})\} \nn \\
&=\sum_{k'=1}^{p^{\circ}}\bigg(
s_{n}^{2}\partial_{\alpha_{k'}^{\circ}}M_{n}(0,\dots,0,\hat{\al}^{\circ}_{n,k'},\dots,\hat{\al}^{\circ}_{n,p^{\circ}},\aes^{\ast},\bes)\hat{\alpha}_{n,k'}^{\circ} \nn\\
&{}\qquad +\frac{s_{n}^{2}}{2}\partial^{2}_{\alpha_{k'}^{\circ}}M_{n}(0,\dots,0,\tilde{\al}^{\circ}_{n,k'},\hat{\al}^{\circ}_{n,k'+1},\dots,\hat{\al}^{\circ}_{n,p^{\circ}},\aes^{\ast},\bes)|\hat{\alpha}_{n,k'}^{\circ}|^{2} + s_{n}^{2}\lam^{a\circ}_{n,k'}R^{a}(\hat{\alpha}_{n,k'}^{\circ})\bigg),
\label{ys:q-nq}
\end{align}
where each point $\tilde{\alpha}_{n,k'}^{\circ}$ is located on the segment connecting zero and $\hat{\alpha}_{n,k'}^{\circ}$. Set $V_{n}(\al)\equiv 0$ and also the rightmost side in \eqref{ys:q-nq} to be $U_{n}(\aes)$. 
By the assumptions, in particular \eqref{hm:for.sparse.LLN}, on an event having probability tending to $1$ we can estimate as
\begin{align}
U_{n}(\aes) &\ge \sum_{k'=1}^{p^{\circ}}\bigg(
\frac{s_{n}^{2}}{2}\partial^{2}_{\alpha_{k'}^{\circ}}M_{n}(0,\dots,0,\tilde{\al}^{\circ}_{n,k'},\hat{\al}^{\circ}_{n,k'+1},\dots,\hat{\al}^{\circ}_{n,p^{\circ}},\aes^{\ast},\bes) -\underline{a}_{k'}
\lam_{\min}\big(\Gam_{0}^{\al^{\circ}}(\al_{0}) \big) \bigg)|\hat{\al}^{\circ}_{n,k'}|^{2} \nn\\
&\ge \sum_{k'=1}^{p^{\circ}}\bigg(
\frac{s_{n}^{2}}{2}\partial^{2}_{\alpha_{k'}^{\circ}}M_{n}(\tz) - \underline{a}_{k'}\lam_{\min}\big(\Gam_{0}^{\al^{\circ}}(\al_{0}) \big) + o_{p}(1) \bigg)|\hat{\al}^{\circ}_{n,k'}|^{2}
\nn\\
&\ge \sum_{k'=1}^{p^{\circ}}\bigg\{\bigg(\frac{1}{2} - \underline{a}_{k'}\bigg) \lam_{\min}\big(\Gam_{0}^{\al^{\circ}}(\al_{0}) \big) + o_{p}(1) \bigg\}|\hat{\al}^{\circ}_{n,k'}|^{2}.
\nn
\end{align}
Hence it follows that for $\underline{a}:=\max_{1\le k'\le p^{\circ}}\underline{a}_{k'}\in(0,1/2)$ we have
\begin{align}
U_{n}(\aes) \ge 
\sum_{k'=1}^{p^{\circ}}\frac{1}{2} \bigg(\frac{1}{2}-\underline{a}_{k'}\bigg) \lam_{\min}\big(\Gam_{0}^{\al^{\circ}}(\al_{0}) \big) |\hat{\al}^{\circ}_{n,k'}|^{2}
\ge \frac{1}{2}\bigg(\frac{1}{2}-\underline{a}\bigg) \lam_{\min}\big(\Gam_{0}^{\al^{\circ}}(\al_{0}) \big) |\aes^{\circ}|^{2},
\nonumber
\end{align}
verifying the second condition in Lemma \ref{Rad08_thm2-} with $\rho=2$: we get $P(\aes^{\circ}=0)\rightarrow 1$. 
To deduce $P(\bes^{\circ}=0)\rightarrow1$, set
\begin{align}
\overline{\mbbh}_{n}(\beta) = t_{n}^{2}\left\{\mbbh_{n}(\aes,\beta)-\mbbh_{n}(\aes,0,\beta^{\ast})\right\}
\nonumber
\end{align}
and follow the same way as above. We have thus obtained the sparse consistency:

\begin{thm}\label{ys:thm_sc}
We have $P(\tes^{\circ}=0)\rightarrow1$ under Assumptions \ref{hm:A_c}, \ref{hm:A_rc} and \ref{ys:sc}.
\end{thm}

%\begin{rem}{\rm 
%Under the tightness of $(\hat{u}_{n},\hat{v}_{n})$ we have
%\begin{align}
%& s_{n}\partial_{\alpha_{k'}^{\circ}}M_{n}(0,\dots,0,\hat{\al}^{\circ}_{n,k'},\dots,\hat{\al}^{\circ}_{n,p^{\circ}},\aes^{\ast},\bes) = O_{p}(1), \qquad k'=1,\dots,p^{\circ},
%\label{ys:a-zero_tight}\\
%& t_{n}\partial_{\beta_{l'}^{\circ}}M_{n}(\aes,0,\dots,0,\hat{\beta}^{\circ}_{n,l'},\dots,\hat{\beta}^{\circ}_{n,q^{\circ}},\hat{\beta}^{\ast}_{n}) = O_{p}(1), \qquad l'=1,\dots,q^{\circ}.
%\label{ys:b-zero_tight}
%\end{align}
%This facts should be effectively used for verification of Assumption \ref{ys:sc}, as in \cite[Section 5]{Rad08}.
%}\end{rem}

\begin{rem}\label{hm:rem_bridge}
\upshape 
Under the tightness of $(\hat{u}_{n},\hat{v}_{n})$ we have
\begin{align}
& s_{n}\partial_{\alpha_{k'}^{\circ}}M_{n}(0,\dots,0,\hat{\al}^{\circ}_{n,k'},\dots,\hat{\al}^{\circ}_{n,p^{\circ}},\aes^{\ast},\bes) = O_{p}(1), \qquad k'=1,\dots,p^{\circ},
\label{ys:a-zero_tight}\\
& t_{n}\partial_{\beta_{l'}^{\circ}}M_{n}(\aes,0,\dots,0,\hat{\beta}^{\circ}_{n,l'},\dots,\hat{\beta}^{\circ}_{n,q^{\circ}},\hat{\beta}^{\ast}_{n}) = O_{p}(1), \qquad l'=1,\dots,q^{\circ}.
\label{ys:b-zero_tight}
\end{align}
This facts should be effectively used for verification of Assumption \ref{ys:sc}, as in \cite[Section 5]{Rad08}. 
As a concrete example, let us consider Assumption \ref{ys:sc} for the bridge-type regularization
\begin{align}
\overline{R}_{n}^{a \circ}(\aes^{\circ})=\sum_{k'=1}^{p^{\circ}}\lambda_{n,k'}^{a\circ}|\hat{\al}^{\circ}_{n,k'}|^{\gamma_{a}}, \quad \overline{R}_{n}^{b \circ}(\bes^{\circ})=\sum_{l'=1}^{q^{\circ}}\lambda_{n,l'}^{b\circ}|\hat{\beta}^{\circ}_{n,l'}|^{\gamma_{b}}, \nn
\end{align}
where $\lambda_{n,k'}^{a\circ},\ \lambda_{n,l'}^{b\circ}>0$ and $\gamma_{a},\ \gamma_{b}\in(0,1]$ are nonrandom tuning parameters with
\begin{align}
\min_{k'}(\lambda_{n,k'}^{a\circ}s_{n}^{\gamma_{a}})\rightarrow\infty,\qquad \min_{l'}(\lambda_{n,l'}^{b\circ}t_{n}^{\gamma_{b}})\rightarrow\infty.
\label{hm:rem.bridge.eq1}
\end{align}
If further $|s_{n}^{-1}\hat{\al}_{n,k'}^{\circ}|=O_{p}(1)$,
% this fails to hold in the Lasso (for linear regression), cf. Zou (2006, Thm.3), hence there is no logical contradiction.
$|t_{n}^{-1}\hat{\beta}_{n,l'}^{\circ}|=O_{p}(1)$, \eqref{ys:a-zero_tight}, and \eqref{ys:b-zero_tight}, then for any $\gamma_{a}\in(0,1]$ we have
\begin{align}
&s_{n}^{2}\partial_{\alpha^{\circ}_{k'}}M_{n}(0,\ldots,0,\hat{\al}^{\circ}_{n,k'},\ldots,\hat{\al}_{n,p^{\circ}}^{\circ},\hat{\al}^{\ast}_{n},\bes)\hat{\al}^{\circ}_{n,k'}+s_{n}^{2}\lambda_{n,k'}^{a\circ}|\hat{\al}^{\circ}_{n,k'}|^{\gamma_{a}} \nn \\
&=\bigg(s_{n}\partial_{\alpha^{\circ}_{k'}}M_{n}(0,\ldots,0,\hat{\al}^{\circ}_{n,k'},\ldots,\hat{\al}_{n,p^{\circ}}^{\circ},\hat{\al}^{\ast}_{n},\bes)\frac{\hat{\al}_{n,k'}^{\circ}}{|\hat{\al}_{n,k'}^{\circ}|}\frac{1}{|s_{n}^{-1}\hat{\al}_{n,k'}^{\circ}|}
\nn\\
&{}\qquad
+s_{n}^{2}|\hat{\al}_{n,k'}^{\circ}|^{-2}\lambda_{n,k'}^{a\circ}|\hat{\al}^{\circ}_{n,k'}|^{\gamma_{a}}\bigg)|\hat{\al}^{\circ}_{n,k'}|^{2} \nn \\
&=\bigg(O_{p}(1)\frac{1}{|s_{n}^{-1}\hat{\al}_{n,k'}^{\circ}|}+\frac{1}{|s_{n}^{-1}\hat{\al}_{n,k'}^{\circ}|^{2}}\lambda_{n,k'}^{a\circ}|\hat{\al}^{\circ}_{n,k'}|^{\gamma_{a}}\bigg)|\hat{\al}_{n,k'}^{\circ}|^{2} \nn \\
&=\bigg(\frac{1}{|s_{n}^{-1}\hat{\al}_{n,k'}^{\circ}|}\left(O_{p}(1)+|s_{n}^{-1}\hat{\al}_{n,k'}^{\circ}|^{\gamma_{a}-1}s_{n}^{\gamma_{a}}\lambda_{n,k'}^{a\circ}\right)\bigg)|\hat{\al}_{n,k'}^{\circ}|^{2}.
\label{hm:rem.bridge.eq2}
\end{align}
Therefore Assumption \ref{ys:sc} holds. For the $\beta$-part, we can follow the same way as above. 
The condition \eqref{hm:rem.bridge.eq1} might turn out to be too much to ask if we look at the coefficient in front of ``$|\hat{\al}_{n,k'}^{\circ}|^{2}$'' in \eqref{hm:rem.bridge.eq2} more precisely, whenever available. 
It should be noted that the above argument may apply for more general regularization terms: see the conditions \eqref{ys:sc-diff-a3} and \eqref{ys:sc-diff-b3} in Section \ref{hm:sec_ergo.diff}.
%as well as Section \ref{hm:sec_tp.selection} for specific choice of regularization terms.
\end{rem}

%%%
\subsubsection{Asymptotic non-degenerate distribution}\label{hm:sec_ad.nz}

We now proceed to derivation of an asymptotic non-trivial distribution of
\begin{align}
\hat{w}_{n}^{\ast}=\left(\hat{u}_{n}^{\ast},\, \hat{v}_{n}^{\ast}\right) := \left( s_{n}^{-1}(\aes^{\ast}-\al^{\ast}_{0}),\, t_{n}^{-1}(\bes^{\ast}-\beta^{\ast}_{0}) \right).
\nonumber
\end{align}
%To this end we will make repeated use of Lemma \ref{Rad08_jal.thm}; although we have already derived the sparse consistency of $\tes^{\circ}$, simultaneous consideration of $\tes^{\circ}$ and $\tes^{\ast}$ via Lemma \ref{Rad08_jal.thm} provides a transparent way to effectively separate the sparse part from the non-sparse one. 
Let Assumptions \ref{hm:A_c}, \ref{hm:A_rc} and \ref{ys:sc} hold, so that we have both $( s_{n}^{-1}(\aes-\al_{0}),\, t_{n}^{-1}(\bes-\beta_{0}))=O_{p}(1)$ and $P(\tes^{\circ}=0)\to 1$ in hand.
%Consider the statistical random field
%\begin{align}
%\mbbm_{n}(u,v) &:= \mbbh_{n}(\al_{0}+s_{n}u,\,\beta_{0}+t_{n}v)-\mbbh_{n}(\theta_{0}) \nn\\
%&= \overline{R}_{n}^{a\circ}(s_{n}v^{\circ}) + \overline{R}_{n}^{b\circ}(t_{n}v^{\circ}) \nn\\
%&{}\qquad +\Big\{M_{n}\left(\tz+(s_{n}u, t_{n}v)\right)-M_{n}(\theta_{0}) \nn \\
%&\qquad\qquad+\big(\overline{R}_{n}^{a\ast}(\al^{\ast}_{0}+s_{n}u^{\ast})-\overline{R}_{n}^{a\ast}(\al_{0}^{\ast})\big)
%+\big(\overline{R}_{n}^{b\ast}(\beta_{0}^{\ast}+t_{n}v^{\ast})-\overline{R}_{n}^{b\ast}(\beta_{0}^{\ast})\big)\Big\}.
%\nonumber
%\end{align}
By the definition, the variable $\hat{w}_{n}^{\ast}$ is a minimizer of the random function
\begin{align}
\mbbm_{n}(u^{\ast},v^{\ast})
&:= \mbbh_{n}(\aes^{\circ},\al_{0}^{\ast}+s_{n}u^{\ast},\bes^{\circ},\beta_{0}^{\ast}+t_{n}v^{\ast})-\mbbh_{n}(\aes^{\circ},\al^{\ast}_{0},\bes^{\circ},\beta^{\ast}_{0}) \nn\\
&= \mbbm_{n}^{0}(u^{\ast},v^{\ast}) + \D\overline{R}_{n}^{a\ast}(u^{\ast}) + \D\overline{R}_{n}^{b\ast}(v^{\ast}),
\nn
\end{align}
where we wrote
\begin{align}
\mbbm_{n}^{0}(u^{\ast},v^{\ast}) &=M_{n}(\aes^{\circ},\al_{0}^{\ast}+s_{n}u^{\ast},\bes^{\circ},\beta_{0}^{\ast}+t_{n}v^{\ast})
-M_{n}(\aes^{\circ},\al^{\ast}_{0},\bes^{\circ},\beta^{\ast}_{0}), \nn\\
\D\overline{R}_{n}^{a\ast}(u^{\ast}) &=\overline{R}_{n}^{a\ast}(\al^{\ast}_{0}+s_{n}u^{\ast})-\overline{R}_{n}^{a\ast}(\al_{0}^{\ast}), \nn\\
\D\overline{R}_{n}^{b\ast}(v^{\ast}) &=\overline{R}_{n}^{b\ast}(\beta_{0}^{\ast}+t_{n}v^{\ast})-\overline{R}_{n}^{b\ast}(\beta_{0}^{\ast}).
\nonumber
\end{align}
Let $D_{n}:=\diag(s_{n}I_{p^{\ast}},t_{n}I_{q^{\ast}})$ and $w^{\ast}:=(u^{\ast},v^{\ast})$. Concerning $\mbbm_{n}^{0}$, we apply the third-order Taylor expansion under the present assumptions to obtain, for a suitable random point $(\tilde{\al}^{\ast}_{n},\tilde{\beta}^{\ast}_{n})$,
\begin{align}
\mbbm_{n}^{0}(u^{\ast},v^{\ast}) &=
D_{n}\p_{\theta^{\ast}}M_{n}(\aes^{\circ},\al_{0}^{\ast},\bes^{\circ},\beta_{0}^{\ast})
[w^{\ast}] + \frac{1}{2}D_{n}\p_{\theta^{\ast}}^{2}M_{n}(\aes^{\circ},\al_{0}^{\ast},\bes^{\circ},\beta_{0}^{\ast})D_{n}[w^{\ast},w^{\ast}]
\nn\\
&{}\qquad + \frac{1}{6}\p_{\theta^{\ast}}^{3}M_{n}(\aes^{\circ},\tilde{\al}^{\ast}_{n},\bes^{\circ},\tilde{\beta}^{\ast}_{n})\left[(D_{n}w^{\ast})^{\otimes 3}\right]
\nn\\
&=\Delta_{n}(\theta_{0})[w^{\ast}] + \frac{1}{2}\Gam_{n}(\tz)[w^{\ast},w^{\ast}] + o_{p}(1),
\label{hm:ast.amn.1}
\end{align}
where
\begin{align}
\Delta_{n}(\theta_{0}) := D_{n}\p_{\theta^{\ast}}M_{n}(\tz), \qquad \Gam_{n}(\tz):=D_{n}\p_{\theta^{\ast}}^{2}M_{n}(\tz)D_{n},
\nonumber
\end{align}
and where the small-order symbol $o_{p}(1)$ in \eqref{hm:ast.amn.1} is valid uniform in $w^{\ast}$ over each compact set in $\mbbr^{p^{\ast}}\times\mbbr^{q^{\ast}}$. Here in the second equality in \eqref{hm:ast.amn.1}, the sparse consistency of $\tes^{\circ}=(\aes^{\circ},\bes^{\circ})$ was used to extract $\D_{n}(\tz)$ and $\Gam_{n}(\tz)$ through the mean-value theorem.

We need the joint weak convergence of these random sequences.

\begin{assum}\label{ys:AL}
There exist random variables $\D_{0}$ and $\Gam_{0}$, and random functions $\D\overline{R}_{0}^{a\ast}(u^{\ast})$ and $\D\overline{R}_{0}^{b\ast}(v^{\ast})$ such that
\begin{align}
\left(\Delta_{n}(\theta_{0}),\, \Gam_{n}(\tz),\, \D\overline{R}_{n}^{a\ast}(\cdot),\, \D\overline{R}_{n}^{b\ast}(\cdot) \right) 
\cil \left( \Delta_{0},\, \Gam_{0},\, \D\overline{R}_{0}^{a\ast}(\cdot),\, \D\overline{R}_{0}^{b\ast}(\cdot) \right)
\nonumber
\end{align}
in $\mathcal{C}(K_{0}\times K_{1})$ for every compact $K_{0}\times K_{1}\subset\mbbr^{p^{\ast}}\times\mbbr^{q^{\ast}}$, 
and that the random function
\begin{align}
\mbbm_{0}(u^{\ast},v^{\ast}):= \Delta_{0}[w^{\ast}] + \frac{1}{2}\Gam_{0}[w^{\ast},w^{\ast}] + \D\overline{R}_{0}^{a\ast}(u^{\ast}) + \D\overline{R}_{0}^{b\ast}(v^{\ast})
\nonumber
\end{align}
has an a.s. unique minimum at $(u^{\ast},v^{\ast})=(\hat{u}_{0}^{\ast},\hat{v}_{0}^{\ast})$.
\end{assum}

The argmax theorem concludes the following.

\begin{thm}\label{hm:thm_AL}
We have $(\hat{u}_{n}^{\ast},\hat{v}_{n}^{\ast})\xrightarrow{\mathcal{L}}(\hat{u}_{0}^{\ast},\hat{v}_{0}^{\ast})$ under Assumptions \ref{hm:A_c}, \ref{hm:A_rc}, \ref{ys:sc} and \ref{ys:AL}.
\end{thm}

Note that under the foregoing assumptions, it follows from \eqref{hm:ast.amn.1} that $\mathbb{M}_{n}^{0}$ has the LAQ structure:
\begin{align}
\mathbb{M}_{n}^{0}(u^{\ast},v^{\ast})=\Delta_{n}(\theta_{0})[w^{\ast}]+\frac{1}{2}\Gamma_{0}[w^{\ast},w^{\ast}]+r_{n}(w^{\ast}),
\nn
%\label{ys:LAQ-3.3.2.}
\end{align}
where the random function $r_{n}$ satisfies that $\sup_{w^{\ast}\in K}|r_{n}(w^{\ast})|\xrightarrow{p}0$ for each compact $K\subset\mathbb{R}^{p^{\ast}+q^{\ast}}$.

\medskip

We now specialize the previous result by focusing on the case of asymptotically mixed-normally distributed $\tes^{\ast}$. To this end, we introduce some further notation and regularity conditions.
For random elements $\mu=\mu(\omega)$ and $\Sig=\Sig(\omega)$ we use the symbol ``$MN(\mu,\Sig)$" for a ``mixed-normal" distribution corresponding to the characteristic function
\begin{align}
z\mapsto E \bigg\{ \exp\bigg( i\mu[u] -\frac{1}{2}\Sig[u,u]\bigg)\bigg\};
\nonumber
\end{align}
we may consider an extension of the original probability space whenever necessary. 
For random variables $X_{n}$ and $X_{0}$, we say that $X_{n}$ ($\mcf$-)stably converges in law to $X_{0}$, denoted by $X_{n}\scl X_{0}$, if $(X_{n},F_{n})\xrightarrow{\mathcal{L}}(X_{0},F_{0})$ for any $\mathcal{F}$-measurable $F_{n}$ such that $F_{n}\cip F_{0}$ (see \cite{JacPro12} for details).

\begin{assum}\label{ys:LAQ_assum}\ 
\begin{enumerate}
\item \label{ys:and-LAQ1}$\Delta_{n}(\theta_{0})\scl \Delta_{0}\sim MN(\mu,\Sigma)$ for some ($\mcf$-measurable) random variables $\mu\in\mathbb{R}^{p^{\ast}+q^{\ast}}$ and $\Sigma\in\mathbb{R}^{p^{\ast}+q^{\ast}}\otimes\mathbb{R}^{p^{\ast}+q^{\ast}}$ with the latter being a.s. positive definite;
\item \label{ys:and-LAQ2}$\Gam_{n}(\tz)\cip\Gam_{0}$ with $P(\Gamma_{0}>0)=1$.
\item \label{ys:and-LAQ3}The functions $R^{a}$ and $R^{b}$ are twice differentiable on $\mathbb{R}\backslash\{0\}$, and for each $k''\in\{1,\ldots,p^{\ast}\}$ and $l''\in\{1,\ldots,q^{\ast}\}$ there correspond nonnegative random variables $\overline{\lambda}_{0,k''}^{a\ast}$ and $\overline{\lambda}_{0,l''}^{b\ast}$ such that $s_{n}\lambda_{n,k''}^{a\ast}\xrightarrow{p}\overline{\lambda}_{0,k''}^{a\ast}$ and $t_{n}\lambda_{n,l''}^{b\ast}\xrightarrow{p}\overline{\lambda}_{0,l''}^{b\ast}$.
\end{enumerate}
\end{assum}

\begin{cor}\label{ys-cor_AL}
Assume that we have Assumptions \ref{hm:A_c}, \ref{hm:A_rc}, \ref{ys:sc}, and \ref{ys:LAQ_assum}. Then, any minimizer $\hat{w}_{n}^{\ast}$ of $\mbbm_{n}$ satisfies that
\begin{align}
\hat{w}_{n}^{\ast} = -\Gamma_{0}^{-1}(\Delta_{n}(\tz)+\Lambda^{\ast}) + o_{p}(1),
\label{hm:linear.rep}
\end{align}
where $\Lambda^{\ast}=(\Lambda^{a\ast},\Lambda^{b\ast})^{\top}$ with
\begin{align}
\Lambda^{a\ast} &:= \Big(\overline{\lambda}_{0,1}^{a\ast}\p R^{a}(\al_{0,1}^{\ast}),\,\ldots\,,\overline{\lambda}_{0,p^{\ast}}^{a\ast}\p R^{a}(\al_{0,p^{\ast}}^{\ast})\Big)^{\top},
\label{hm:Lam.a_def} \\
\Lambda^{b\ast} &:= \Big(\overline{\lambda}_{0,1}^{b\ast}\p R^{b}(\beta_{0,1}^{\ast}),\,\ldots\,,\overline{\lambda}_{0,q^{\ast}}^{b\ast}\p R^{b}(\beta_{0,q^{\ast}}^{\ast})\Big)^{\top}.
\label{hm:Lam.b_def}
\end{align}
In particular, we have the asymptotic mixed normality:
\begin{align}
(\hat{u}_{n}^{\ast},\hat{v}_{n}^{\ast}) \cil (\hat{u}_{0}^{\ast},\hat{v}_{0}^{\ast}) \sim MN\left(-\Gamma_{0}^{-1}(\mu+\Lambda^{\ast}),\Gamma_{0}^{-1}\Sigma\Gamma_{0}^{-1}\right).
\nonumber
\end{align}
\end{cor}

\begin{proof}
Write $F_{n}(\theta^{\ast})= \mbbh_{n}(\aes^{\circ},\al^{\ast},\bes^{\circ},\beta^{\ast}) -  \mbbh_{n}(\aes^{\circ},\al^{\ast}_{0},\bes^{\circ},\beta^{\ast}_{0})$. Under the assumptions, on an event whose probability tending to $1$ it follows from the mean-value theorem that
\begin{align}
\left( -D_{n}\p_{\theta^{\ast}}^{2}F_{n}(\check{\theta}_{n})D_{n}\right) [\hat{w}_{n}^{\ast}] = D_{n}\p_{\theta^{\ast}}F_{n}(\tz^{\ast})
\label{hm:amn.add1}
\end{align}
for some random point $\check{\theta}_{n}$ such that $\check{\theta}_{n}\cip\tz^{\ast}$. Moreover, for suitable points $(\check{\al}^{\ast}_{n},\check{\beta}^{\ast}_{n})$,
\begin{align}
D_{n}\p_{\theta^{\ast}}F_{n}(\tz^{\ast}) &= \D_{n}(\tz) + \Lam^{\ast} + o_{p}(1),
\nn\\
-D_{n}\p_{\theta^{\ast}}^{2}F_{n}(\check{\theta}_{n})D_{n} &= -\Gam_{n}(\tz) - \diag\left(
s_{n}^{2}\p_{\al^{\ast}}^{2}\overline{R}_{n}^{a \ast}(\check{\al}^{\ast}_{n}),\, t_{n}^{2}\p_{\beta^{\ast}}^{2}\overline{R}_{n}^{b \ast}(\check{\beta}^{\ast}_{n})
\right) + o_{p}(1) \nn\\
&= -\Gam_{0} + o_{p}(1),
\nonumber
\end{align}
where, as we did in \eqref{hm:ast.amn.1} before, the sparse consistency of $\tes^{\circ}$ was used to extract $\D_{n}(\tz)$ and $\Gam_{n}(\tz)$ in the right-hand sides. 
Hence \eqref{hm:linear.rep} follows from \eqref{hm:amn.add1}. Then the latter claim is trivial.
\end{proof}

%%%
\subsubsection{Uniform tail-probability estimate}\label{unif}

We now proceed to the tail-probability estimate of the scaled estimator $\hat{w}_{n}:=(\hat{u}_{n},\hat{v}_{n})$, the main concern of this paper; recall the definition \eqref{hm:uhat.vhat.def}.
To this end, we will make use of the PLDI result in \cite{Yos11}, which is, though not specified, originally developed for the LAQ statistical random fields. Throughout this section we assume that the matrix $\Gam_{0}$ is a.s. constant.
%we focus on ergodic-type models where the matrix $\Gam_{0}$ is a.s. constant.
Write $u=(u^{\circ},u^{\ast})$ and $v=(v^{\circ},v^{\ast})$.
We will proceed with the multi-step PLDI argument in \cite{Yos11} (see also \cite[Section 4]{Mas_LM}), which roughly goes as follows:

\begin{itemize}
\item At the first step we consider the mixed-rates random field 
\begin{align}
\mbbm_{n}(u,\beta)&:=\mathbb{H}_{n}\left(\alpha_{0}+s_{n}u,\beta\right)-\mathbb{H}_{n}(\alpha_{0},\beta) \nonumber \\
&=M_{n}\left(\alpha_{0}+s_{n}u,\beta\right)-M_{n}(\alpha_{0},\beta)+\overline{R}_{n}^{a\circ}(s_{n}u^{\circ})+\overline{R}_{n}^{a\ast}(\alpha_{0}^{\ast}+s_{n}u^{\ast})-\overline{R}_{n}^{a\ast}(\alpha_{0}^{\ast}) \nonumber \\
&=:s_{n}^{-2}f_{n}(u^{\circ})+t_{n}^{-2}g_{n}(u,\beta),
\nn
\end{align}
where
\begin{align}
s_{n}^{-2}f_{n}(u^{\circ})&:=\overline{R}_{n}^{a\circ}(s_{n}u^{\circ}), \nonumber \\
t_{n}^{-2}g_{n}(u,\beta)&:=M_{n}\left(\alpha_{0}+s_{n}u,\beta\right)-M_{n}(\alpha_{0},\beta)+\overline{R}_{n}^{a\ast}(\alpha_{0}^{\ast}+s_{n}u^{\ast})-\overline{R}_{n}^{a\ast}(\alpha_{0}^{\ast}) \nonumber
\end{align}
with regarding $\beta$ as a nuisance parameter. We will derive the $\alpha$-PLDI
\begin{align}
\sup_{r}\sup_{n}r^{L}P(|\hat{u}_{n}|\geq r)<\infty
\nonumber
\end{align}
through proving
\begin{align}
\sup_{r>0}\sup_{n>0}r^{L}P\left[\sup_{(u,\beta)\in\{|u|\geq r\}\times\Theta_{\beta}}\{-\mbbm_{n}(u,\beta)\}\geq 0\right]<\infty
\label{ys:a-PLDI}
\end{align}
for any $L>0$, according to the facts that the probability in \eqref{ys:a-PLDI} is equal to or greater than $P(|\hat{u}_{n}|\geq r)$, that $\mbbm_{n}(0,\beta)=0$, and that $\sup_{\beta\in\Theta_{\beta}}\{\mbbh_{n}(\al_{0},\beta)-\mbbh_{n}(\aes,\beta)\}\geq0$. 

\item At the second step, by plugging-in $\hat{\alpha}_{n}$ we look at the mixed-rates random field
\begin{align}
\mbbm_{n}(v)&:=\mathbb{H}_{n}\left(\hat{\alpha}_{n},\beta_{0}+t_{n}v\right)-\mathbb{H}_{n}(\hat{\alpha}_{n},\beta_{0}) \nonumber \\
&=M_{n}\left(\hat{\alpha}_{n},\beta_{0}+t_{n}v\right)-M_{n}(\hat{\alpha}_{n},\beta_{0})+\overline{R}_{n}^{b\circ}(t_{n}v^{\circ})+\overline{R}_{n}^{b\ast}(\beta_{0}^{\ast}+t_{n}v^{\ast})-\overline{R}_{n}^{b\ast}(\beta_{0}^{\ast}) \nonumber \\
&=:s_{n}^{-2}f_{n}(v^{\circ})+t_{n}^{-2}g_{n}(v),
\nonumber
\end{align}
where
\begin{align}
s_{n}^{-2}f_{n}(v^{\circ})&:=\overline{R}_{n}^{b\circ}(t_{n}v^{\circ}), \nonumber \\
t_{n}^{-2}g_{n}(v)&:=M_{n}\left(\hat{\alpha}_{n},\beta_{0}+t_{n}v\right)-M_{n}(\hat{\alpha}_{n},\beta_{0})+\overline{R}_{n}^{b\ast}(\beta_{0}^{\ast}+t_{n}v^{\ast})-\overline{R}_{n}^{b\ast}(\beta_{0}^{\ast}). \nonumber 
\end{align}
Then, we will derive the $\beta$-PLDI
\begin{align}
\sup_{r}\sup_{n}r^{L}P(|\hat{v}_{n}|\geq r)<\infty
\nonumber
\end{align}
from
\begin{align}
\sup_{r>0}\sup_{n>0}r^{L}P\left[\sup_{|v|\geq r}\{-\mbbm_{n}(v)\}\geq 0\right]<\infty.
\label{ys:b-PLDI}
\end{align}
Combining \eqref{ys:a-PLDI} and \eqref{ys:b-PLDI} yields the joint-PLDI
\begin{align}
\sup_{r>0}\sup_{n>0}r^{L}P(|\hat{w}_{n}|\geq r)<\infty,
\label{ys:joint-PLDI}
\end{align}
so that the uniform tail-probability estimates of $\hat{w}_{n}$ follows.
\end{itemize}

We see that the verification of $\alpha$-PLDI reduces to that of the PLDI of the form
\begin{align}
\sup_{r>0}\sup_{n>0}r^{L}P\left[\sup_{(u,\beta)\in\{|u|\ge r\}\times\Theta_{\beta}}\{-t_{n}^{-2}g_{n}(u,\beta)\}\ge 0\right]<\infty,
\label{ys:al2-PLDI}
\end{align}
since we have 
\begin{align}
\inf_{u\in\{|u|\ge r\}}\overline{R}_{n}^{a\circ}(s_{n}u^{\circ})\geq0 \quad {\rm a.s.}
\nn
\end{align}
and therefore
\begin{align}
\sup_{(u,\beta)\in\{|u|\ge r\}\times\Theta_{\beta}}\{-\mathbb{M}_{n}(u,\beta)\}\ge 0 &\iff \sup_{(u,\beta)\in\{|u|\ge r\}\times\Theta_{\beta}}\{-t_{n}^{-2}g_{n}(u,\beta)\}\ge \inf_{u\in\{|u|\ge r\}}s_{n}^{-2}f_{n}(u^{\circ})
\nonumber \\
&\ \Longrightarrow \sup_{(u,\beta)\in\{|u|\ge r\}\times\Theta_{\beta}}\{-t_{n}^{-2}g_{n}(u,\beta)\}\ge 0.
%\label{ys:a-est}
\nn
\end{align}
Typically $-t_{n}^{-2}g_{n}(u,\beta)$ fails to be LAQ because of the presence of the regularization term $\overline{R}_{n}^{a\ast}(\alpha_{0}^{\ast}+s_{n}u^{\ast})-\overline{R}_{n}^{a\ast}(\alpha_{0}^{\ast})$.
However, we can apply the argument in \cite{Yos11} with regarding $-t_{n}^{-2}g_{n}(u,\beta)$ as ``LAQ-like" structure (see \eqref{ys:LAQ_Ad}--\eqref{ys:LAQ_Ar} below). 
Similarly, the PLDI 
\begin{align}
\sup_{r>0}\sup_{n>0}r^{L}P\left[\sup_{v\in\{|v|\geq r\}}\{-t_{n}^{2}g_{n}(v)\}\geq0\right]<\infty 
\label{ys:b-PLDI2}
\end{align}
verifies the $\beta$-PLDI through \eqref{ys:b-PLDI} since we have
\begin{align}
\inf_{v\in\{|v|\ge r\}}\overline{R}_{n}^{b\circ}(t_{n}v^{\circ})\geq0 \quad {\rm a.s.} \nn
\end{align}

\medskip

First we set some regularity conditions on the part $M_{n}$. 

\begin{assum}\label{pldi:ass1}\ 
\begin{enumerate}
\item There exist nonrandom functions $\widetilde{M}_{0}^{a}: \Theta_{\alpha}\rightarrow\mathbb{R}$ and $\widetilde{M}_{0}^{b}: \Theta_{\beta}\rightarrow\mathbb{R}$, and positive constants $\delta_{1}^{a}$, $\delta_{1}^{b}$, $\chi^{a}$ and $\chi^{b}$ such that for all $K>0$,
\begin{itemize}
\item $\sup_{n>0}E\left[\sup_{\theta\in\Theta}\left|s_{n}^{-2\delta_{1}^{a}}\left\{s_{n}^{2}\big(M_{n}(\alpha,\beta)-M_{n}(\alpha_{0},\beta)\big)-\widetilde{M}_{0}^{a}(\alpha)\right\}\right|^{K}\right]<\infty$.
\item $\sup_{n>0}E\left[\sup_{\beta\in\Theta_{\beta}}\left|t_{n}^{-2\delta_{1}^{b}}\left\{t_{n}^{2}\big(M_{n}(\al_{0},\beta)-M_{n}(\al_{0},\beta_{0})\big)-\widetilde{M}_{0}^{b}(\beta)\right\}\right|^{K}\right]<\infty$.
\item $\widetilde{M}_{0}^{a}(\alpha)\geq \chi^{a}|\alpha-\alpha_{0}|^{2}$, \quad $\widetilde{M}_{0}^{b}(\beta)\geq \chi^{b}|\beta-\beta_{0}|^{2}$.
\end{itemize}
\item There exist nonrandom matrices $C_{0}(\beta)$ and $C_{0}>0$, and constants $\delta_{2}^{a},\ \delta_{2}^{b}\in(0,1/2]$ such that for all $K>0$,
\begin{itemize}
\item $\inf_{\beta}\lambda_{\min}(C_{0}(\beta))>0$.
\item $\sup_{n>0}E\left[\sup_{\beta\in \Theta_{\beta}}\left(s_{n}^{-2\delta_{2}^{a}}\left|s_{n}^{2}\partial_{\alpha}^{2}M_{n}\left(\alpha_{0},\beta\right)-C_{0}(\beta)\right|\right)^{K}\right]<\infty$.
\item $\sup_{n>0}E\left[\left(t_{n}^{-2\delta_{2}^{b}}\left|t_{n}^{2}\partial_{\beta}^{2}M_{n}\left(\al_{0},\beta_{0}\right)-C_{0}\right|\right)^{K}\right]<\infty$.
\item $\sup_{n>0}E\left[\sup_{\al\in\Theta_{\al}}\left|s_{n}t_{n}^{2(1-\delta_{2}^{b})}\p_{\al}\p_{\beta}^{2}M_{n}(\al,\beta_{0})\right|^{K}\right]<\infty$.
\end{itemize}
\item For all $K>0$,
\begin{itemize}
\item $\sup_{n>0}E\left[\sup_{\beta\in\Theta_{\beta}}\left|s_{n}\partial_{\alpha}M_{n}(\alpha_{0},\beta)\right|^{K}\right]<\infty$, \quad $\sup_{n>0}E\left[\left|t_{n}\partial_{\beta}M_{n}(\al_{0},\beta_{0})\right|^{K}\right]<\infty$.
\item $\sup_{n>0}E\left[\sup_{\theta\in\Theta}\left|s_{n}^{2}\partial_{\alpha}^{3}M_{n}(\theta)\right|^{K}\right]<\infty$, \quad $\sup_{n>0}E\left[\sup_{\theta\in\Theta}\left|t_{n}^{2}\partial_{\beta}^{3}M_{n}(\theta)\right|^{K}\right]<\infty$.
\item $\sup_{n>0}E\left[\sup_{\al\in\Theta_{\al}}|s_{n}t_{n}\p_{\al}\p_{\beta}M_{n}(\al,\beta_{0})|^{K}\right]<\infty$.
\end{itemize}
\end{enumerate}
\end{assum}

%\medskip
%{\color{red}
%$t_{n}^{2}\p_{\beta}^{2}M_{n}(\aes,\beta_{0})+C_{0}=t_{n}^{2}\p_{\beta}^{2}M_{n}(\al_{0},\beta_{0})+C_{0}+t_{n}^{2}\p_{\al}\p_{\beta}^{2}M_{n}(\tilde{\al},\beta_{0})[\aes-\al_{0}]$\\
%$\therefore\ \sup_{n>0}E[\sup_{\al\in\Theta_{\al}}|s_{n}t_{n}^{2(1-\delta_{2}^{b})}\p_{\al}\p_{\beta}^{2}M_{n}(\al,\beta_{0})|^{K}]<\infty.$  e.g. diffusion: $s_{n}t_{n}^{2(1-\delta)}=(nh_{n})^{-1}n^{\delta-1/2}h_{n}^{\delta}$.}
%{\color{red}
%\begin{enumerate}
%\item $t_{n}\partial_{\beta}M_{n}(\hat{\alpha}_{n},\beta_{0})=t_{n}\partial_{\beta}M_{n}(\al_{0},\beta_{0})+t_{n}\p_{\al}\partial_{\beta}M_{n}(\tilde{\al},\beta_{0})[\aes-\al_{0}]$\\
%$\therefore\ \sup_{n>0}E[\sup_{\al\in\Theta_{\al}}|s_{n}t_{n}\p_{\al}\p_{\beta}M_{n}(\al,\beta_{0})|^{K}]<\infty.$  e.g. diffusion: $s_{t}t_{n}=(nh_{n})^{-1}h_{n}^{1/2}$.
%\end{enumerate}
%}

Assumption \ref{pldi:ass1} is borrowed from \cite[Theorem 3(c)]{Yos11}. 
%Our primary stance here is to give a set of conditions 
Although the description may seem to be stringent at first glance, quite often we may verify them in a straightforward manner under suitable stability conditions such as the (exponential) ergodicity and the boundedness of moments; see \cite[Section 6]{Yos11} and also Section \ref{hm:sec_ergo.diff} below.
Further, we should note that the present regularity conditions could be relaxed in compensation for more messy description;
see Remark \ref{hm:rem_laq.nondble}.

\medskip

We also need moment conditions on the regularization terms.

\begin{assum}\label{pldi:ass7}
There exist constants $\nu^{a},\ \nu^{b}\in(0,1/2)$ such that the following conditions hold for any $K>0$,
\begin{enumerate}
\item $\ds{\sup_{n>0}E\left[\sup_{\alpha\in\Theta_{\al}}\left(s_{n}^{1+2\nu^{a}}\overline{R}_{n}^{a}(\al)\right)^{K}\right]<\infty}$, \quad $\ds{\sup_{n>0}E\left[\sup_{\beta\in\Theta_{\beta}}\left(t_{n}^{1+2\nu^{b}}\overline{R}_{n}^{b}(\beta)\right)^{K}\right]<\infty}$;
\item $\ds{\sup_{n>0}\max_{k''}E[|s_{n}\lambda_{n,k''}^{a\ast}|^{K}]<\infty}$, \quad $\ds{\sup_{n>0}\max_{l''}E[|t_{n}\lambda_{n,l''}^{b\ast}|^{K}]<\infty}$. 
\end{enumerate}
\end{assum}

The following theorem gives us the joint-PLDI \eqref{ys:joint-PLDI}.

\begin{thm}\label{thm:pldi}
For any $L>0$, \eqref{ys:joint-PLDI} holds under Assumptions \ref{pldi:ass1} and \ref{pldi:ass7}. Additionally if we have the weak convergence 
$\left(\hat{u}^{\circ}_{n},\hat{u}_{n}^{\ast},\hat{v}^{\circ}_{n},\hat{v}_{n}^{\ast}\right) \xrightarrow{\mathcal{L}} \left(\hat{u}_{0}^{\circ},\hat{u}_{0}^{\ast},\hat{v}_{0}^{\circ},\hat{v}_{0}^{\ast} \right)$ for some random vector $\left(\hat{u}_{0}^{\circ},\hat{u}_{0}^{\ast},\hat{v}_{0}^{\circ},\hat{v}_{0}^{\ast} \right)$, then the moment convergence
\begin{align}
E[f(\hat{u}_{n},\hat{v}_{n})]\rightarrow E[f(\hat{u}_{0}^{\circ},\hat{u}_{0}^{\ast},\hat{v}_{0}^{\circ},\hat{v}_{0}^{\ast})]
\label{eq:p.g.}
\end{align}
holds for all continuous $f:\mbbr^{p+q}\to\mbbr$ of at most polynomial growth. 
\end{thm}

\begin{proof}
At the first step, we will derive the PLDI \eqref{ys:al2-PLDI}.
The Taylor expansion around $u=0$ yields 
\begin{align}
-t_{n}^{-2}g_{n}(u,\beta) = \Delta_{n}(\beta)[u] - \frac{1}{2}C_{0}(\beta)[u,u] + r_{n}(u,\beta),
\nonumber
\end{align}
where
%
%\begin{align}
%t_{n}^{-2}g_{n}(u,\beta)&=s_{n}\partial_{\alpha}M_{n}(\alpha_{0},\beta)\left[u\right]+\frac{1}{2}C_{0}(\beta)[u,u] \nonumber \\
%&\quad +\int_{0}^{1}(1-s)\left\{s_{n}^{2}\partial_{\alpha}^{2}M_{n}\left(\alpha_{0}+s_{n}su,\beta\right)+C_{0}(\beta)\right\}ds\left[u,u\right] \nn \\
%&\quad +\{\overline{R}_{n}^{a\ast}(\alpha_{0}^{\ast}+s_{n}u^{\ast})-\overline{R}_{n}^{a\ast}(\alpha_{0}^{\ast})\}, \nonumber
%\end{align}
%so that for the $\alpha$-PLDI we have the LAQ-like expression
%
%\begin{align}
%\mathbb{M}_{n}(u,\beta;\theta_{0})&=\mathbb{H}_{n}(\al_{0}+A_{n}^{11}(\theta_{0})u,\beta)-\mathbb{H}_{n}(\al_{0},\beta) \nn \\
%&=\Delta_{n}(\beta;\theta_{0})[u]+\frac{1}{2}\Gamma_{0}(\beta;\theta_{0})[u,u]+r_{n}(u,\beta;\theta_{0}) \quad u\in\mathbb{R}^{p},
%\nn%\label{eq:plaq}
%\end{align}
%
%
%...
%
%where
\begin{align}
\Delta_{n}(\beta)&:=-s_{n}\partial_{\alpha}M_{n}(\alpha_{0},\beta), \label{ys:LAQ_Ad}\\
%\Gamma_{0}(\beta)&=C_{0}(\beta), \label{ys:LAQ_Ag}\\
r_{n}(u,\beta)&:=-\int_{0}^{1}(1-s)\left\{s_{n}^{2}\partial_{\alpha}^{2}M_{n}\left(\alpha_{0}+s_{n}su,\beta\right)-C_{0}(\beta)\right\}ds\left[u,u\right] \nn\\
&{}\qquad -\{\overline{R}_{n}^{a\ast}(\alpha_{0}^{\ast}+s_{n}u^{\ast})-\overline{R}_{n}^{a\ast}(\alpha_{0}^{\ast})\}. 
\label{ys:LAQ_Ar}
\end{align}
We will check the conditions [A1], [A4], [A6], [B1] and [B2] of \cite[Theorem 1]{Yos11} (see Appendix \ref{appendix-A} for details).

From Assumptions \ref{pldi:ass1} and \ref{pldi:ass7}, [B1] is trivial and for any $\theta\in\Theta$ we have
\begin{align}
\mathbb{Y}_{n}^{a}(\alpha,\beta)&:=-s_{n}^{2}\{\mathbb{H}_{n}(\alpha,\beta)-\mathbb{H}_{n}(\alpha_{0},\beta)\} \nonumber \\
&=-s_{n}^{2}\{M_{n}(\alpha,\beta)-M_{n}(\alpha_{0},\beta)\}-s_{n}^{2}\{\overline{R}_{n}^{a}(\al)-\overline{R}_{n}^{a\ast}(\al_{0}^{\ast})\} \nonumber \\
&\xrightarrow{p}-\widetilde{M}_{0}^{a}(\alpha)=:\mathbb{Y}_{0}^{a}(\alpha)\leq -\chi^{a}|\alpha-\alpha_{0}|^{2},
\nn
\end{align}
hence [B2] also holds. We will derive [A1] and [A6], along with controlling the auxiliary parameters given in [A4]; for clarity, we will put the superscript ``$a$'' or ``$b$'' for them.
In the sequel, we will write $x_{n}\lesssim y_{n}$ if $\sup_{n}(x_{n}/y_{n})<\infty$.

First, we will verify [A1]. We have
\begin{align}
&\int_{0}^{1}(1-s)\left\{s_{n}^{2}\partial_{\alpha}^{2}M_{n}\left(\alpha_{0}+s_{n}su,\beta\right)-C_{0}(\beta)\right\}ds\left[u,u\right] \nonumber \\
&=\int_{0}^{1}(1-s)\left\{s_{n}^{2}\partial_{\alpha}^{2}M_{n}\left(\alpha_{0},\beta\right)-C_{0}(\beta)+\int_{0}^{1}s_{n}^{2}\partial_{\alpha}^{3}M_{n}\left(\alpha_{0}+s_{n}stu,\beta\right)dt\left[s_{n}su\right]\right\}ds\left[u,u\right], \nonumber
\end{align}
hence
\begin{align}
\frac{|r_{n}(u,\beta)|}{1+|u|^{2}}&\lesssim \frac{|u|^{2}}{1+|u|^{2}}\left|s_{n}^{2}\partial_{\alpha}^{2}M_{n}\left(\alpha_{0},\beta\right)-C_{0}(\beta)\right|
\nn\\
&\quad+\frac{|u|^{2}}{1+|u|^{2}}s_{n}|u|\int_{0}^{1}\int_{0}^{1}\left|s_{n}^{2}\partial_{\alpha}^{3}M_{n}\left(\alpha_{0}+s_{n}stu,\beta\right)\right|dtds \nonumber \\
&\quad+\frac{1}{1+|u|^{2}}\left|\overline{R}_{n}^{a\ast}(\alpha_{0}^{\ast}+s_{n}u^{\ast})-\overline{R}_{n}^{a\ast}(\alpha_{0}^{\ast})\right|.
\label{eq:pldi1}
\end{align} 
Let
\begin{align}
U_{n}^{a}(r) &:= \big\{ u\in\mbbr^{p};\, \alpha_{0}+s_{n}u\in\Theta_{\alpha},\quad r\leq|u|\leq s_{n}^{-(1-\zeta^{a})} \big\},
%U_{n}^{b}(r) &= \big\{ v\in\mbbr^{q};\, \beta_{0}+t_{n}v\in\Theta_{\beta},\quad r\leq|v|\leq t_{n}^{-(1-\zeta^{b})}\big\}
\nonumber
\end{align}
where $\zeta^{a}\in(0,1)$.
From \eqref{eq:pldi1} and the Chebyshev's inequality, the probability in [A1] can be estimated as follows:
\begin{align}
&P\left(\sup_{(u,\beta)\in U_{n}^{a}(r)\times\Theta_{\beta}}\frac{|r_{n}(u,\beta)|}{1+|u|^{2}}\geq r^{-\rho_{1}^{a}}\right) \nonumber \\
&\lesssim r^{d\rho_{1}^{a}}\left\{E\left[\sup_{(u,\beta)\in U_{n}^{a}(r)\times\Theta_{\beta}}\left(\frac{|u|^{2}}{1+|u|^{2}}\left|s_{n}^{2}\partial_{\alpha}^{2}M_{n}\left(\alpha_{0},\beta\right)-C_{0}(\beta)\right|\right)^{d}\right]\right. \nonumber \\
&\quad+E\left[\sup_{(u,\beta)\in U_{n}^{a}(r)\times\Theta_{\beta}}\left(\frac{|u|^{2}}{1+|u|^{2}}s_{n}|u|\int_{0}^{1}\int_{0}^{1}\left|s_{n}^{2}\partial_{\alpha}^{3}M_{n}\left(\alpha_{0}+s_{n}stu,\beta\right)\right|dtds\right)^{d}\right] \nonumber \\
&\quad\left.+E\left[\sup_{u\in U_{n}^{a}(r)}\left(\frac{1}{1+|u|^{2}}\left|\overline{R}_{n}^{a\ast}(\alpha_{0}^{\ast}+s_{n}u^{\ast})-\overline{R}_{n}^{a\ast}(\alpha_{0}^{\ast})\right|\right)^{d}\right]\right\}
\label{eq:pldi2}
\end{align}
for some $d>0$. Letting $\zeta^{a}\in(0,\min\{\delta_{1}^{a},2\delta_{2}^{a},1/2\})$ and $\xi\in(0,\zeta^{a}/(1-\zeta^{a}))$. Note that $-\zeta^{a}/2+(1-\zeta^{a})\xi/2<0$, and therefore $-\delta_{2}^{a}+(1-\zeta^{a})\xi/2<0$. 
Then, for the first term of the right-hand side in \eqref{eq:pldi2}, it follows from Assumption \ref{pldi:ass1} that
\begin{align}
&\sup_{n>0}E\left[\sup_{(u,\beta)\in U_{n}^{a}(r)\times\Theta_{\beta}}\left(\frac{|u|^{2}}{1+|u|^{2}}\left|s_{n}^{2}\partial_{\alpha}^{2}M_{n}\left(\alpha_{0},\beta\right)-C_{0}(\beta)\right|\right)^{d}\right] \nonumber \\
&=\sup_{n>0}\left\{E\left[\sup_{\beta\in \Theta_{\beta}}\left(s_{n}^{-2\delta_{2}^{a}}\left|s_{n}^{2}\partial_{\alpha}^{2}M_{n}\left(\alpha_{0},\beta\right)-C_{0}(\beta)\right|\right)^{d}\right]\sup_{u\in U_{n}^{a}(r)}\left(\frac{|u|^{2}}{1+|u|^{2}}s_{n}^{2\delta_{2}^{a}}|u|^{\xi}|u|^{-\xi}\right)^{d}\right\} \nonumber \\
&\lesssim r^{-d\xi}\sup_{n>0}s_{n}^{-2\{-\delta_{2}^{a}+(1-\zeta^{a})\xi/2\}d}\lesssim r^{-d\xi}.
\label{eq:pldi3}
\end{align}
Further, for the second term of the right-hand side in \eqref{eq:pldi2}, under Assumption \ref{pldi:ass1} we obtain the estimate:
\begin{align}
&\sup_{n>0}E\left[\sup_{(u,\beta)\in U_{n}^{a}(r)\times\Theta_{\beta}}\left(\frac{|u|^{2}}{1+|u|^{2}}s_{n}|u|\int_{0}^{1}\int_{0}^{1}\left|s_{n}^{2}\partial_{\alpha}^{3}M_{n}\left(\alpha_{0}+s_{n}stu,\beta\right)\right|dtds\right)^{d}\right] \nonumber \\
&\lesssim \sup_{n>0}\sup_{u\in U_{n}^{a}(r)}\left(\frac{|u|^{2}}{1+|u|^{2}}s_{n}|u||u|^{\xi}|u|^{-\xi}\right)^{d}\lesssim r^{-d\xi}\sup_{n>0}s_{n}^{-2\{-\zeta^{a}/2+(1-\zeta^{a})\xi/2\}d}\lesssim r^{-d\xi}.
\label{eq:pldi4}
\end{align}
As for the third term of the right-hand side in \eqref{eq:pldi2}, the local Lipschitz continuity \eqref{hm:llc} and Assumption \ref{pldi:ass7} yield
\begin{align}
& \sup_{n>0}E\left[\sup_{u\in U_{n}^{a}(r)}\left(\frac{1}{1+|u|^{2}}\left|\overline{R}_{n}^{a\ast}(\alpha_{0}^{\ast}+s_{n}u^{\ast})-\overline{R}_{n}^{a\ast}(\alpha_{0}^{\ast})\right|\right)^{d}\right]
\nn\\
&{}\qquad\lesssim\sup_{n>0}E\left[\sup_{u\in U_{n}^{a}(r)}\left(\frac{1}{1+|u|^{2}}\sum_{k''=1}^{p^{\ast}}\lambda_{n,k''}^{a\ast}s_{n}|u^{\ast}_{k''}|\right)^{d}\right] \nonumber \\
&{}\qquad\lesssim \sup_{n>0}E\left[\sup_{u\in U_{n}^{a}(r)}\sum_{k''=1}^{p^{\ast}}\left(\frac{\lambda_{n,k''}^{a\ast}s_{n}|u^{\ast}_{k''}|}{1+|u|^{2}}\right)^{d}\right]\lesssim r^{-d}.
\label{eq:pldi5}
\end{align}
Fix a $\rho_{1}^{a}\in(0,\xi)$. Then from \eqref{eq:pldi2}--\eqref{eq:pldi5}, we get for any $L>0$
\begin{align}
\sup_{r>0}\sup_{n>0}r^{L}P\left(\sup_{(u,\beta)\in U_{n}^{a}(r)\times\Theta_{\beta}}\frac{|r_{n}(u,\beta)|}{1+|u|^{2}}\geq r^{-\rho_{1}^{a}}\right)<\infty.
\nn
\end{align}

For condition [A4], we take $\rho_{2}^{a}>2\zeta^{a}$ and $\nu_{1}^{a}>0$ such that $1-2\nu_{1}^{a}-\rho_{2}^{a}>0$. Then, the condition [A6] is trivial since we have Assumptions \ref{pldi:ass1} and \ref{pldi:ass7}, and
\begin{align}
\mathbb{Y}_{n}^{a}(\alpha,\beta)-\mathbb{Y}_{0}^{a}(\alpha)=-s_{n}^{2}\{M_{n}(\alpha,\beta)-M_{n}(\alpha_{0},\beta)\}-s_{n}^{2}\{\overline{R}_{n}^{a}(\alpha)-\overline{R}_{n}^{a\ast}(\alpha_{0}^{\ast})\}+\widetilde{M}_{0}^{a}(\alpha),
\nonumber
\end{align}
i.e. we get 
\begin{align}
& \sup_{n>0}E\left[\left(\sup_{\beta\in\Theta_{\beta}}|\Delta_{n}(\beta)|\right)^{N_{1}}\right]<\infty, \nn\\
& \sup_{n>0}E\left[\left(\sup_{\theta\in\Theta}s_{n}^{-2(1/2-\nu_{1}^{a})}|\mathbb{Y}_{n}^{a}(\theta)-\mathbb{Y}_{0}^{a}(\alpha)|\right)^{N_{2}}\right]<\infty,
\nonumber
\end{align}
for $N_{1}=L(1-\rho_{1}^{a})^{-1}>0$ and $N_{2}=L(1-2\nu_{1}^{a}-\rho_{2}^{a})^{-1}>0$. Therefore [A4] holds with taking the tuning  parameters in [A4] as above. The proof of the $\alpha$-PLDI \eqref{ys:al2-PLDI} is complete.

\medskip

At the second step of the multistep-PLDI argument for showing the $\beta$-PLDI \eqref{ys:b-PLDI2}, we can take an analogous way to the derivation of the $\al$-PLDI. We have
\begin{align}
-t_{n}^{-2}g_{n}(v) = \D_{n}[u] -\frac{1}{2}C_{0}[u,u] + r_{n}(v),
\nonumber
\end{align}
where
\begin{align}
\Delta_{n}&:=-t_{n}\partial_{\beta}M_{n}(\hat{\alpha}_{n},\beta_{0}), \nn \\
%\Gamma_{0}&=C_{0}, \nn \\
r_{n}(v)&:=-\int_{0}^{1}(1-s)\left\{t_{n}^{2}\partial_{\beta}^{2}M_{n}\left(\hat{\alpha}_{n},\beta_{0}+t_{n}sv\right)-C_{0}\right\}ds\left[v,v\right]-\{\overline{R}_{n}^{b\ast}(\beta_{0}^{\ast}+t_{n}v^{\ast})-\overline{R}_{n}^{b\ast}(\beta_{0}^{\ast})\}. \nonumber
\end{align}
Verification of [B1] and [B2] are the same as before: for any $\beta\in\Theta_{\beta}$, we have
\begin{align}
\mathbb{Y}_{n}^{b}(\beta)&:=-t_{n}^{2}\{\mbbh_{n}(\hat{\alpha}_{n},\beta)-\mbbh_{n}(\hat{\alpha}_{n},\beta_{0})\} \nonumber \\
&=-t_{n}^{2}\{M_{n}(\hat{\alpha}_{n},\beta)-M_{n}(\hat{\alpha}_{n},\beta_{0})\}-t_{n}^{2}\{\overline{R}^{b}_{n}(\beta)-\overline{R}_{n}^{b\ast}(\beta_{0}^{\ast})\} \nonumber \\
&\xrightarrow{p}-\widetilde{M}_{0}^{b}(\beta)=:\mathbb{Y}_{0}^{b}(\beta)\leq -\chi^{b}|\beta-\beta_{0}|^{2}.
\nn
\end{align}
Through expanding the $ds$-integrand of $r_{n}(v)$, we get the following as in \eqref{eq:pldi1}:
%\begin{align}
%&\int_{0}^{1}(1-s)\left\{t_{n}^{2}\partial_{\beta}^{2}M_{n}\left(\hat{\alpha}_{n},\beta_{0}+t_{n}sv\right)+C_{0}\right\}ds\left[v,v\right] \nonumber \\
%&=\int_{0}^{1}(1-s)\left\{t_{n}^{2}\partial_{\beta}^{2}M_{n}\left(\hat{\alpha}_{n},\beta_{0}\right)+C_{0}+\int_{0}^{1}t_{n}^{2}\partial_{\beta}^{3}M_{n}\left(\hat{\alpha}_{n},\beta_{0}+t_{n}stv\right)dt\left[t_{n}sv\right]\right\}ds\left[v,v\right], \nonumber
%\end{align}
%hence we obtain
\begin{align}
\frac{|r_{n}(v)|}{1+|v|^{2}}&\lesssim \frac{|v|^{2}}{1+|v|^{2}}\left|t_{n}^{2}\partial_{\beta}^{2}M_{n}\left(\hat{\alpha}_{n},\beta_{0}\right)-C_{0}\right|+\frac{|v|^{2}}{1+|v|^{2}}t_{n}|v|\int_{0}^{1}\int_{0}^{1}\left|t_{n}^{2}\partial_{\beta}^{3}M_{n}\left(\hat{\alpha}_{n},\beta_{0}+t_{n}stv\right)\right|dtds \nonumber \\
&\quad+\frac{1}{1+|v|^{2}}\left|\overline{R}_{n}^{b\ast}\left(\beta_{0}+t_{n}v^{\ast}\right)-\overline{R}_{n}^{b\ast}\left(\beta_{0}^{\ast}\right)\right|. 
\nn%\label{eq:b-pldi1}
\end{align}
Introducing $U_{n}^{b}(r) = \big\{ v\in\mbbr^{q};\, \beta_{0}+t_{n}v\in\Theta_{\beta},~ r\leq|v|\leq t_{n}^{-(1-\zeta^{b})}\big\}$ for $\zeta^{b}\in(0,1)$, we can follow the same line as in the derivation of \eqref{eq:pldi2}--\eqref{eq:pldi5} as before, to verify that for any $L>0$
\begin{align}
\sup_{r>0}\sup_{n>0}r^{L}P\left(\sup_{v\in U_{n}^{b}(r)}\frac{|r_{n}(v)|}{1+|v|^{2}}\geq r^{-\rho_{1}^{b}}\right)<\infty,
\nn
\end{align}
for some number $\rho_{1}^{b}\in(0,\xi)$ with $\xi\in(0,\zeta^{b}/(1-\zeta^{b}))$ for $\zeta^{b}\in(0,\min\{\delta_{1}^{b},2\delta_{2}^{b},1/2\})$. 
Note that we made use of the moment bound
\begin{align}
\sup_{n>0}E\left[\left|s_{n}^{-1}(\aes-\al_{0})\right|^{K}\right]<\infty\ {\rm for\ all}\ K>0,
\nonumber
\end{align}
ensured by the inequality $\sup_{r}\sup_{n}r^{L}P(|\hat{u}_{n}|\geq r)<\infty$ which has been previously obtained from \eqref{ys:al2-PLDI}. 
Verification of [A6] along with [A4] can be done in the same ways as before. The proof of the $\beta$-PLDI \eqref{ys:b-PLDI2} is complete, concluding \eqref{ys:joint-PLDI}. The latter claim of the theorem is trivial.
\end{proof}

\begin{rem}{\rm 
We note that, in view of \cite[Theorem 1]{Yos11}, the twice differentiability of the LAQ part $M_{n}$ and the constancy of $\Gam_{0}$ are not essential. All the assertions presented in this section can also go for possibly non-differential $M_{n}$ as long as statistical random fields associated with $M_{n}$ is of LAQ: $M_{n}(\tz + A_{n}u) - M_{n}(\tz) = \D_n [u] +\frac{1}{2}\Gam_{0}[u,u] + r_{n}(u)$ with (quasi-)score sequence $\D_{n}$, asymptotic (quasi-)information matrix $\Gam_{0}$, and remainder term $r_{n}(u)=o_{p}(1)$ (locally uniformly in $u$).
The resulting set of conditions becomes somewhat less concise.
}\label{hm:rem_laq.nondble}
\end{rem}

%%%%%
\subsection{Standard case}\label{hm:sec_st.asymp.}
In Section \ref{hm:sec_s.asymp} we have focused on the sparse asymptotics, where the underlying statistical random field is of mixed-rates type. Nevertheless, it is possible to treat the standard asymptotics as well, where the localization via matrix \eqref{hm:intro1}:
\begin{align}
\mathbb{M}_{n}(u;\theta_{0})=\mathbb{H}_{n}(\theta_{0}+A_{n}(\theta_{0})u)-\mathbb{H}_{n}(\theta_{0}),
\nn
\end{align}
can completely determine the asymptotic distribution of $\tes$. In this case, we cannot take the sparse consistency into account even if the true parameter $\tz$ has zero components. We note that the random field $\mathbb{M}_{n}(\cdot;\theta_{0})$ still may not be LAQ according to choice of the regularization terms; for example, the non-sparse asymptotics for the linear regression considered in \cite[Theorem 2]{KniFu00} belongs to this standard-asymptotics category, while, for example, the ridge-type regularization \cite{HoeKen70} under suitable stochastic-order conditions on the regularization coefficients $\lam^{a}_{n,k}$ and $\lam^{b}_{n,l}$ (recall \eqref{hm:regu.term.add1}) leads to the LAQ structure (most often the locally asymptotically mixed-normal one).

In the standard asymptotics, the underlying statistical random field is normally not of the mixed-rates type. It is of course possible and easier to state results analogous to Theorems \ref{hm:thm_consistency}, \ref{hm:thm_rate}, \ref{hm:thm_AL}, Corollary \ref{ys-cor_AL}, and Theorem \ref{thm:pldi}, simply by removing the factor $\theta^{\circ}$, hence in particular Assumption \ref{ys:sc}, in the foregoing arguments.

Let us be more specific. We consider \eqref{hm:Hn} and \eqref{hm:regu.term.add1} for the random function $\mbbh_{n}$ consisting of $M_{n}(\theta)$, $\overline{R}^{a}_{n}(\al)$, and $\overline{R}^{b}_{n}(\beta)$ as before, except that we here do {\it not} impose the factorization $\theta=(\theta^{\circ},\theta^{\ast})$; the true value $\tz$ may be such that $\theta_{0,k}\ne 0$ for all $k\in\{1,\dots,p+q\}$, and we no longer consider the decompositions \eqref{hm:Ra_form} and \eqref{hm:Rb_form}. 
Then, by inspecting the corresponding proofs in Section \ref{hm:sec_s.asymp}, it is straightforward to observe the following.
\begin{itemize}
\item \textit{The assertion of Theorem \ref{hm:thm_consistency} holds as it is.}
\item Define Assumption \ref{hm:A_rc}${}^{\star}$ to be the same one as Assumption \ref{hm:A_rc} except that Assumption \ref{hm:A_rc}.2 therein is replaced by ``$s_{n}\lam^{a}_{n,k}=O_{p}(1)$ and $t_{n}\lam^{b}_{n,l}=O_{p}(1)$ for each $k$ and $l$''. 
Then, \textit{the assertion of Theorem \ref{hm:thm_rate} holds under Assumptions \ref{hm:A_c} and \ref{hm:A_rc}${}^{\star}$.}
\item Turning to asymptotic distribution of $(\hat{u}_{n}, \hat{v}_{n}) = ( s_{n}^{-1}(\aes-\al_{0}), t_{n}^{-1}(\bes-\beta_{0}) )$, we first note that under Assumptions \ref{hm:A_c} and \ref{hm:A_rc}${}^{\star}$ we have
\begin{align}
\mbbm_{n}(u,v) &:=\mbbh_{n}(\al_{0}+s_{n}u,\beta_{0}+t_{n}v)-\mbbh_{n}(\al_{0},\beta_{0}) \nn\\
&=\Delta_{n}(\theta_{0})[w] + \frac{1}{2}\Gam_{n}(\tz)[w,w] + \D\overline{R}_{n}^{a}(u) + \D\overline{R}_{n}^{b}(v) + o_{p}(1),
\nn
\end{align}
where $w:=(u,v)$ and we wrote, different from Section \ref{hm:sec_ad.nz}, $\Delta_{n}(\theta_{0}) := D_{n}\p_{\theta}M_{n}(\tz)$ and $\Gam_{n}(\tz):=D_{n}\p_{\theta}^{2}M_{n}(\tz)D_{n}$ for $D_{n}:=\diag(s_{n}I_{p},\,t_{n}I_{q})$, and also $\D\overline{R}_{n}^{a}(u):=\overline{R}_{n}^{a}(\al_{0}+s_{n}u)-\overline{R}_{n}^{a}(\al_{0})$ and $\D\overline{R}_{n}^{b}(v):=\overline{R}_{n}^{b}(\beta_{0}+t_{n}v)-\overline{R}_{n}^{b}(\beta_{0})$. Let us define Assumption \ref{ys:AL}${}^{\star}$ to be the same one as Assumption \ref{ys:AL} with all the symbol ``$\ast$'' therein removed from the notation. Then, as a variant of Theorem \ref{hm:thm_AL}, we have the following claim:
\textit{we have $(\hat{u}_{n},\hat{v}_{n})\xrightarrow{\mathcal{L}}(\hat{u}_{0},\hat{v}_{0})$ under Assumptions \ref{hm:A_c}, \ref{hm:A_rc}${}^{\star}$, and \ref{ys:AL}${}^{\star}$.} Of course, the asymptotic distribution may not be mixed-normal, as in the non-Gaussian limit specified in \cite{KniFu00} for the linear regression model. An analogue to Corollary \ref{ys-cor_AL} can be deduced similarly under appropriate conditions.

\item Finally as for Theorem \ref{thm:pldi}, we partly modify Assumption \ref{pldi:ass7} by removing the symbol ``$\ast$'' therein from the notation, and call it Assumption \ref{pldi:ass7}${}^{\star}$. Then the claim is as follows:
\textit{under Assumptions \ref{pldi:ass1} and \ref{pldi:ass7}${}^{\star}$ we have \eqref{ys:joint-PLDI} for any $L>0$, and additionally if we have the weak convergence $(\hat{u}_{n},\hat{v}_{n}) \xrightarrow{\mathcal{L}} (\hat{u}_{0},\hat{v}_{0})$ for some random vector $(\hat{u}_{0},\hat{v}_{0})$, then the moment convergence
\begin{align}
E[f(\hat{u}_{n},\hat{v}_{n})]\rightarrow E[f(\hat{u}_{0},\hat{v}_{0})]
\nn
\end{align}
holds for all continuous $f:\mbbr^{p+q}\to\mbbr$ of at most polynomial growth.}

\end{itemize}

%%%%%
\subsection{Discussion}

\subsubsection{Specific choice of regularization term}\label{hm:sec_tp.selection}

Let us recall \eqref{hm:Ra_form} and \eqref{hm:Rb_form}: $\overline{R}_{n}^{a}(\al) =\overline{R}_{n}^{a \circ}(\al^{\circ}) + \overline{R}_{n}^{a \ast}(\al^{\ast}) = \sum_{k'=1}^{p^{\circ}}\lam^{a \circ}_{n,k'}R^{a}(\al^{\circ}_{k'}) + \sum_{k''=1}^{p^{\ast}}\lam^{a \ast}_{n,k''}R^{a}(\al^{\ast}_{k''})$ and $\overline{R}_{n}^{b}(\beta) =\overline{R}_{n}^{b \circ}(\beta^{\circ}) + \overline{R}_{n}^{b \ast}(\beta^{\ast}) = \sum_{l'=1}^{q^{\circ}}\lam^{b \circ}_{n,l'}R^{b}(\beta^{\circ}_{l'}) + \sum_{l''=1}^{q^{\ast}}\lam^{b \ast}_{n,l''}R^{b}(\beta^{\ast}_{l''})$. 
As a simple explicit component-wise regularization, we remark on a formal selection of the regularization terms $\overline{R}^{a}_{n}(\al)$ and $\overline{R}^{b}_{n}(\beta)$ in terms of maximizing the posterior probability in Bayesian setting. 
We will think of the LAQ part $M_{n}(\theta)$ as the minus of a (quasi) log-likelihood quantity, with the prior distribution of $(\al,\beta)$ of mutually independent components. Namely, we consider prior density, say $\pi^{a}(\al;\lam^{a})$ and $\pi^{b}(\beta;\lam^{b})$ of $\al$ and $\beta$, respectively, of the forms
\begin{align}
\pi^{a}(\al;\lam^{a}) &= \prod_{k'=1}^{p^{\circ}}C^{a}(\lam^{a\circ}_{k'})\exp\big(  -\lam^{a\circ}_{k'}R^{a}(\al_{k'}^{\circ})\big) \cdot 
\prod_{k''=1}^{p^{\ast}}C^{a}(\lam^{a\ast}_{k''})\exp\big(  -\lam^{a\ast}_{k''}R^{a}(\al_{k''}^{\ast})\big), \nn\\
\pi^{b}(\beta;\lam^{b}) &= \prod_{l'=1}^{q^{\circ}}C^{b}(\lam^{b\circ}_{l'})\exp\big(  -\lam^{b\circ}_{l'}R^{b}(\beta_{l'}^{\circ})\big) \cdot 
\prod_{l''=1}^{q^{\ast}}C^{b}(\lam^{b\ast}_{l''})\exp\big(  -\lam^{b\ast}_{l''}R^{b}(\beta_{l''}^{\ast})\big).
\nonumber
\end{align}
We then set $\overline{R}^{a}_{n}(\al)$ and $\overline{R}^{b}_{n}(\beta)$ to be 
$-\log \pi^{a}(\al;\lam_{n}^{a})$ and $-\log \pi^{b}(\beta;\lam_{n}^{b})$ up to additive constants which do not depend on $(\al,\beta)$. 
At this stage, the functions $R^{a}$ and $R^{b}$, and the regularization intensity (possibly random) sequences 
$\lam_{n}^{a}=\{(\lam^{a \circ}_{n,k'})_{k'},\, (\lam^{a \ast}_{n,k''})_{k''}\}$ and $\lam_{n}^{b}=\{(\lam^{b \circ}_{n,l'})_{l'},\, (\lam^{b \ast}_{n,l''})_{l''}\}$
are left unspecified. 
The functions $R^{a}$ and $R^{b}$ determine the class of regularization (such as the ridge type, the lasso type, and so on). 
We suppose they are given and discuss about choice of $\lam_{n}^{a}$ and $\lam_{n}^{b}$ in computing $\tes=\tes(\lam)$; 
how to choose them is in general an unavoidable and very often annoying problem in regularized estimation, 
and there are a lot of possibilities of reasoning on this point from different points of view.\footnote{Historically, development of statistical model-shrinkage methodologies goes back to the so-called $L^{0}$-type regularization. Most popular among others are the classical model-selection criteria Akaike's AIC (see Section \ref{hm:sec_aic}) and Schwarz's BIC, both of which provide theoretical bases of how to choose $\lam$ in significant-parameter selection, or how to fix a significance level in the context of testing hypothesis in a nested model. These came from different contexts, although the random functions to be optimized take rather similar forms: ``the principal part (typically log-likelihood)'' plus ``a regularization term''. See e.g. \cite{BurAnd02} for details.
} 

As in \cite{WanLiGuo07}, we take the root of the following equation with respect to $\lam=(\lam^{a},\lam^{b})$:
\begin{align}
& \p_{\lam}\mbbh_{n}(\theta) = \left(-\p_{\lam^{a}}\log\pi^{a}(\al;\lam^{a}),\, -\p_{\lam^{b}}\log\pi^{b}(\beta;\lam^{b})\right)=0 \nn\\
&\iff \p_{\lam^{a\circ}_{k'}}C^{a}(\lam^{a\circ}_{k'}) = R^{a}(\al_{k'}^{\circ}),\quad 
\p_{\lam^{a\ast}_{k''}}C^{a}(\lam^{a\ast}_{k''}) = R^{a}(\al_{k''}^{\ast}),\nn\\
& \phantom{\iff}\quad \p_{\lam^{b\circ}_{l'}}C^{b}(\lam^{b\circ}_{l'}) = R^{b}(\beta_{l'}^{\circ}),\quad 
\p_{\lam^{b\ast}_{l''}}C^{b}(\lam^{b\ast}_{l''}) = R^{b}(\beta_{l''}^{\ast}).
\label{hm:lam_post}
\end{align}
That is to say, we formally regard $\lam=\lam(\theta)$ maximizing the posterior as an optimal choice. 
We do not know how to set the value $\theta$ \textit{a priori}, and want to use the non-regularized estimator $\tilde{\theta}_{n}=(\tilde{\theta}_{n,k})_{k}:=\tes(0)$; 
then we plug-in the estimate $\tilde{\theta}_{n}$ into \eqref{hm:lam_post} and solve the $p+q$ simultaneous equations with respect to $\lam$, which produces a simple explicit data-driven choice of $\lam_{n}$.

\medskip

Let us take a closer look at the exponential power distribution, which is quite popular in regularized estimation and corresponds to the bridge-type regularization:
\begin{align}
\pi(\theta;\lam)=\prod_{k=1}^{p}\frac{\gam\lam_{k}^{1/\gam}}{2\Gam(1/\gam)}\exp(-\lam_{k}|\theta_{k}|^{\gam})
\nonumber
\end{align}
for some (given) $\gam>0$, where $\lam=(\lam_{k})_{k\le p}$, $\lam_{k}\ge 0$. We here only deal with $\gamma\leq1$. In this case, the identity $\frac{\p}{\p\lam_{k}}\log\pi(\theta;\lam)=0$ leads to the choice
\begin{align}
\lam_{k}=\frac{1}{\gam|\theta_{k}|^{\gam}}.
\label{hm:lnc.1}
\end{align}
Note that this simple choice may fail to meet our technical conditions in the following two points; let us only look at the case of $\al_{k}$, for the case of $\beta_{l}$ can be argued similarly. 
\begin{itemize}
\item The first point is the stochastic-order conditions on the regularization terms. Suppose that an $s_{n}^{-1}$-consistent initial estimator $(\tilde{\al}_{n,k})_{k}$ is available, then \eqref{hm:lnc.1} results in\footnote{For a random variable $X_{n}$, we write $X_{n}\xrightarrow{p}+\infty$ if $P(X_{n}>K)\rightarrow 1$ for every $K>0$.}
\begin{align}
{\lam}_{n,k}^{a}=\frac{1}{\gam|\tilde{\al}_{n,k}|^{\gam}}\left\{
\begin{array}{ll}
\text{$=O_{p}(1)$ and not $o_{p}(1)$} & \text{if $\al_{0,k}\ne 0$,} \\[2mm]
=O_{p}(s_{n}^{-\gam}) & \text{if $\al_{0,k}=0$ and $s_{n}^{-1}\tilde{\al}_{n,k}\ne o_{p}(1)$,} \\[2mm]
\cip +\infty & \text{if $\al_{0,k}=0$ and $s_{n}^{-1}\tilde{\al}_{n,k}= o_{p}(1)$,}
\end{array}\right.
\label{hm:sa.lam}
\end{align}
for $k=1,\dots,p$; typically, the first two cases are available. We should note that the rate of convergence $s_{n}^{-1}$ of $(\tilde{\al}_{n,k})_{k}$ is not essential: a different rate will just change divergence rate of $\lam_{n,k}^{a}$ for $\al_{0,k}=0$ in \eqref{hm:sa.lam}.
\item The second one is that reciprocals of $|\tilde{\al}_{n,k}|$ may not have finite high-order moments, which does matter when trying to verify Assumption \ref{pldi:ass7}: indeed, the typical case, where $\sqrt{n}(\tilde{\al}_{n,k}-\al_{0,k})$ is asymptotically normally distributed, only moments of order less than $1$ exist. 
\end{itemize}
Below we discuss about these two points in more details.

\medskip

Concerning the {\it first point}, let us first note that the first two cases in \eqref{hm:sa.lam} are comparable to the standard asymptotics (e.g. \cite{KniFu00}), since then the random field $(u,v)\mapsto \mbbh_{n}(\al_{0}+s_{n}u,\beta_{0}+t_{n}v)-\mbbh_{n}(\al_{0},\beta_{0})$ can completely describe the asymptotic distribution of the regularized estimator $(\aes,\bes)$.
To consider the sparse asymptotics, we need to modify \eqref{hm:lnc.1} more or less artificially, such as
\begin{align}
\lam_{n,k}^{a}=\frac{l_{n}}{\gam|\tilde{\al}_{n,k}|^{\gam}}
\label{hm:sa.lam.2}
\end{align}
for some sequence $l_{n}\to\infty$ (such as $\log n$, as seen below), as was already mentioned \cite{WanLiGuo07}. We do not have any justification for this ad-hoc choice \eqref{hm:sa.lam.2}, however, it may enable us to apply, for example, the bridge type regularization even when we do not know beforehand if the true value $\al_{0}$ indeed has a zero component or not, since each regularization weight $\lam_{n,k}^{a}$ asymptotically tells us if $\al_{0,k}$ is zero or not. To show this, assume additionally that
\begin{align}
s_{n}^{\gam}l_{n} \to 0.
\nonumber
\end{align}
Then, for \eqref{hm:sa.lam.2} we have
\begin{align}
s_{n}^{\gam}\lam_{n,k}^{a}\left\{
\begin{array}{ll}
=O_{p}(s_{n}^{\gam}l_{n})\cip 0 & \text{if $\al_{0,k}\ne 0$,} \\[2mm]
\ds{=\frac{l_{n}}{\gam|s_{n}^{-1}\tilde{\al}_{n,k}|^{\gam}}\cip +\infty} & \text{if $\al_{0,k}= 0$.}
\end{array}
\right.
\nonumber
\end{align}
Note that the condition \eqref{hm:rem.bridge.eq1} in Remark \ref{hm:rem_bridge} is fulfilled, and also the condition $s_{n}\lam_{n,k}^{a}=o_{p}(1)$ for $\al_{0,k}\ne 0$ as well. 
%Of course, this reasoning is asymptotic and at this stage we still do not have any clear view about finite-sample performance. 
Thus, the modification \eqref{hm:sa.lam.2} suffices if we are only concerned with the weak convergence plus the sparse consistency of $\aes$, without the tail-probability estimate.

From the above discussion, we can conclude the following claim concerning the conditions on the regularization terms:
\begin{cor}
Let $\tilde{\al}_{n}=(\tilde{\al}_{n,k})_{k}$ and $\tilde{\beta}_{n}=(\tilde{\beta}_{n,l})_{l}$ be $s_{n}^{-1}$-consistent and $t_{n}^{-1}$-consistent estimators of $\al$ and $\beta$, respectively,
such that $s_{n}^{-1}\tilde{\al}_{n,k}\ne o_{p}(1)$ and that $t_{n}^{-1}\tilde{\beta}_{n,l}\ne o_{p}(1)$. 
If we set \eqref{hm:regu.term.add1} with 
\begin{align}
&R^{a}(\al_{k})=|\al_{k}|^{\gamma^{a}}, \quad R^{b}(\beta_{l})=|\beta_{l}|^{\gamma^{b}}, 
\label{ys:bridge1}\\
&\lambda_{n,k}^{a}=\frac{l_{n}^{a}}{\gamma^{a}|\tilde{\al}_{n,k}|^{\gamma^{a}}},\quad \lambda_{n,l}^{b}=\frac{l_{n}^{b}}{\gamma^{b}|\tilde{\beta}_{n,l}|^{\gamma^{b}}}
\nn
\end{align}
for some $l_{n}^{a},\ l_{n}^{b}\to\infty$ and $\gamma^{a},\ \gamma^{b}\in(0,1]$ such that
\begin{align}
s_{n}^{\gam^{a}}l_{n}^{a} \to 0,\quad t_{n}^{\gam^{b}}l_{n}^{b}\to0,
\label{ys:bridge2}
\end{align}
then Assumptions \ref{hm:A_c}.3, \ref{hm:A_rc}.2 and \ref{ys:LAQ_assum}.3 hold. If further we have \eqref{ys:a-zero_tight} and \eqref{ys:b-zero_tight}, then Assumption \ref{ys:sc} holds.
\end{cor}

\medskip

To further remedy the {\it second point}, where Assumption \ref{pldi:ass7} comes into play in, we may for example consider a variant of \eqref{hm:sa.lam.2} of the form
\begin{align}
\lam_{n,k}^{a}=\frac{l_{n}}{\gam(|\tilde{\al}_{n,k}| \vee \ep^{a}_{n,k})^{\gam}}
\nonumber
\end{align}
for a positive nonrandom sequence $\ep_{n,k}^{a}$ such that
\begin{align}
\ep^{a}_{n,k} \asymp s_{n}
\nonumber
\end{align}
for $k=1,\dots,p$, in the sense that the sequences $(\ep^{a}_{n,k}/s_{n})$ and $(s_{n}/\ep^{a}_{n,k})$ are bounded. Then
\begin{align}
|\lam_{n,k}^{a}|\lesssim l_{n}I(\al_{0,k}\ne 0) + s_{n}^{-\gam}l_{n} I(\al_{0,k}=0) \lesssim s_{n}^{-\gam}l_{n},
\nonumber
\end{align}
so that Assumption \ref{pldi:ass7} is satisfied if
\begin{align}
\text{$\gam\in(0,1)$ \quad and \quad $l_{n}s_{n}^{1-\gam}\lesssim 1$.}
\nonumber
\end{align}
Indeed, we then have $\sup_{\al_{k}}\left|s_{n}^{1+2\nu_{a}}\lam_{n,k}^{a}R^{a}(\al_{k})\right| + |s_{n}\lam_{n,k}^{a}| 
\lesssim l_{n}s_{n}^{1+2\nu_{a}-\gam} + l_{n}s_{n}^{1-\gam} \lesssim l_{n}s_{n}^{1-\gam}\lesssim 1$.

A possible practical choice of $\ep^{a}_{n,k}$ would be
\begin{align}
\ep^{a}_{n,k}=C_{a}s_{n}
\nonumber
\end{align}
for a fixed constant $C_{a}>0$, which can be arbitrary but would more or less affect finite-sample performance. The asymptotic theory does not tell us how to choose $C_{a}$, remaining as a non-trivial matter. Nevertheless, in practice we may for example set $C_{a}$ so small that $C_{a}s_{n}$ equals a threshold, say $10^{-5}$, for which we basically have to make a subjective choice. Also, we may regard an estimate of $\al_{k}$ as exactly $0$ if its estimate is less than $\ep^{a}_{n,k}$. Although this kind of thresholding seems to be quite artificial, it would be necessary in numerical optimization, which does not in principle return exact zero but just only very small number even when the target value equals zero.

From the above discussion about the moment conditions on the regularization terms, we conclude the following claim:
\begin{cor}
Let $\tilde{\al}_{n}=(\tilde{\al}_{n,k})_{k}$ and $\tilde{\beta}_{n}=(\tilde{\beta}_{n,l})_{l}$ be $s_{n}^{-1}$-consistent and $t_{n}^{-1}$-consistent estimators of $\al$ and $\beta$, respectively,
such that $s_{n}^{-1}\tilde{\al}_{n,k}\ne o_{p}(1)$ and that $t_{n}^{-1}\tilde{\beta}_{n,l}\ne o_{p}(1)$.
If we set \eqref{hm:regu.term.add1} with \eqref{ys:bridge1} and
\begin{align}
\lambda_{n,k}^{a}=\frac{l_{n}^{a}}{\gamma^{a}(|\tilde{\al}_{n,k}|\vee \epsilon_{n,k}^{a})^{\gamma^{a}}},\quad \lambda_{n,l}^{b}=\frac{l_{n}^{b}}{\gamma^{b}(|\tilde{\beta}_{n,l}|\vee\epsilon_{n,l}^{b})^{\gamma^{b}}},
\nn
\end{align}
where $\epsilon_{n,k}^{a}$ and $\epsilon_{n,l}^{b}$ are positive nonrandom sequences such that
\begin{align}
\ep^{a}_{n,k} \asymp s_{n}, \quad \ep^{b}_{n,l} \asymp t_{n}
\nonumber
\end{align}
for any $k\in\{1,\ldots,p\}$ and $l\in\{1,\ldots,q\}$ and where $l_{n}^{a},\ l_{n}^{b}\to\infty$ and $\gamma^{a},\ \gamma^{b}\in(0,1)$ satisfy \eqref{ys:bridge2} and
\begin{align}
s_{n}^{1-\gam^{a}}l_{n}^{a} \to 0,\quad t_{n}^{1-\gam^{b}}l_{n}^{b}\to 0,
\nonumber
\end{align}
then Assumptions \ref{hm:A_c}.3, \ref{hm:A_rc}.2, \ref{ys:LAQ_assum}.3 and \ref{pldi:ass7} hold. If further we have \eqref{ys:a-zero_tight} and \eqref{ys:b-zero_tight}, Assumption \ref{ys:sc} holds.
\end{cor}

\medskip

As has been well recognized, non-convex optimization such as the case of the bridge regularization with $\gam<1$ can be a hard task. Toward stabilization of numerical computation, we could incorporate existing general devices for manipulating $R^{a}$, such as the local quadratic approximation \cite{FanLi01}.

%%%
\subsubsection{Prediction-related issues}
\label{hm:sec_aic}

Based on our moment-convergence results, it is straightforward to verify many kinds of prediction-error asymptotics, such as the asymptotic behavior of the mean-squared error $\E(|\tes-\tz|^{2})=V_{0}/n+o(1/n)$ in its typical form; among many others, we refer to \cite{ChaIng11} and \cite{IngYan14} for some related rigorous treatments.
We may also apply our moment-convergence results to validate the celebrated AIC methodology (\cite{Aka74}, and also \cite{FinWei02}, \cite{UchYos01}, and \cite{UchYos06}) even under the sparse asymptotics. Recent studies exactly in this direction contain \cite{Shi16} and \cite{UmeShiMasNin15}, where the uniform integrability of the sparse maximum-likelihood estimator with the bridge-like regularization played an important role for validating the asymptotic bias correction. Below we will briefly discuss how we can extend the result of \cite{UmeShiMasNin15} to cover a broader range of statistical models with locally asymptotically normal structure. To keep things simple, we are only concerned here with correctly specified models,
keeping the multi-scaling setup with $\theta=(\al,\beta)$; again we note that cases of single scaling can be handled as well without any essential change, just by ignoring the $\beta$-part.

\medskip

We keep adopting the basic setup and notation introduced in Section \ref{hm:sec_setup}. In what follows, we focus on a regularized minus log-likelihood function \eqref{hm:Hn}:
\begin{align}
\mathbb{H}_{n}(\theta)=\mathbb{H}_{n}(\alpha,\beta)=M_{n}(\al,\beta)+\overline{R}_{n}^{a}(\al) +\overline{R}_{n}^{b}(\beta),
\nonumber
\end{align}
where $\{M_{n}(\theta)\}_{\theta\in\Theta}$ is a {\it correctly specified minus log-likelihood function}, namely, $M_{n}(\tz)$ corresponds to the minus likelihood function of true data-generating model.
Recall the notation: $\D_{n}(\tz)=D_{n}\p_{\theta^{\ast}}M_{n}(\tz)$ and $\Gam_{n}(\theta)=D_{n}\p_{\theta^{\ast}}^{2}M_{n}(\theta)D_{n}$.
Building on the results in Section \ref{hm:sec_s.asymp}, we here presuppose that $M_{n}$ is smooth enough in $\theta\in\Theta$, and also that the following conditions hold.
\begin{itemize}
\item $\D_{n}(\tz)\cil N_{p^{\ast}+q^{\ast}}(0,\Gam_{0})$ with the asymptotic covariance matrix $\Gam_{0}$ being nonrandom and positive-definite, $\Gam_{n}(\tz)\cip\Gam_{0}$, and both $\{\D_{n}(\tz)\}_{n}$ and $\{\sup_{\theta}|\Gam_{n}(\theta)|\}_{n}$ are $L^{r}(P)$-bounded for sufficiently large $r>0$.

\item We have the sparse consistency for the zero-parameters: $P(\tes^{\circ}=0)\to 1$.
\item Asymptotic normality for the non-zero parameters:
\begin{align}
\hat{w}_{n}^{\ast}:=\left(s_{n}^{-1}(\aes^{\ast}-\al_{0}^{\ast}),\, t_{n}^{-1}(\bes^{\ast}-\beta_{0}^{\ast}) \right)
=-\Gamma_{0}^{-1}(\Delta_{n}(\tz)+\Lambda^{\ast}) + o_{p}(1),
\label{hm:linear.rep_2}
\end{align}
where $\Lam^{\ast}$ is a constant vector.
\item The uniform probability estimate
\begin{align}
\sup_{r>0}\sup_{n}r^{L}P\left(|\hat{w}_{n}|\ge r\right)<\infty
\nonumber
\end{align}
for some $L>2$ is valid for $\hat{w}_{n}:=(s_{n}^{-1}(\aes-\al_{0}),\, t_{n}^{-1}(\bes-\beta_{0}))$.
\end{itemize}
The right hand side of \eqref{hm:linear.rep_2} tends in distribution to the normal distribution $N_{p^{\ast}+q^{\ast}}(-\Gam_{0}^{-1}\Lam^{\ast}, \Gam_{0}^{-1})$;
apart from the factor $\Lam^{\ast}$, this is the typical form of asymptotically efficient maximum-likelihood estimation, where $\Gam_{0}$ is the Fisher information matrix of the non-zero parameter part. 

Let us denote by $X^{n}$ data for constructing the estimator $\tes=\tes(X^{n})$, and $\tilde{X}^{n}$ an independent copy of $X^{n}$. For clarity, we will often write $M_{n}(\theta)=M_{n}(\theta;X^{n})$. The expected Kullback-Leibler divergence from the true data-generating distribution to the prediction model for $\tilde{X}^{n}$ based on the regularized estimator $\tes$ is given by
\begin{align}
E\left[\tilde{E}\left\{-M_{n}(\tz;\tilde{X}^{n})\right\} - \tilde{E}\left\{-M_{n}(\tes(X^{n});\tilde{X}^{n})\right\}\right],
\nonumber
\end{align}
where $\tilde{E}$ denotes the expectation with respect to $\tilde{X}^{n}$. In principle, it is desirable to make this quantity, equivalently $E\otimes\tilde{E}\{M_{n}(\tes(X^{n});\tilde{X}^{n})\}$,
as small as possible. A natural estimator for the last quantity is $M_{n}(\tes(X^{n});X^{n})$, but it is well-known that the statistics has a significant bias due to the double usage of original data $X^{n}$ for both estimation and prediction. The bias corrected version necessitates evaluation of the quantity
\begin{align}
M_{n}(\tes(X^{n});X^{n}) - E\left[M_{n}(\tes(X^{n});X^{n})-\tilde{E}\left\{M_{n}(\tes(X^{n});\tilde{X}^{n})\right\}\right],
\nonumber
\end{align}
but unfortunately, it is not a statistics in general. The crucial step in the classical AIC methodology for relative model assessment in terms of the best Kullback-Leibler divergence based prediction is to derive and/or approximate the bias term
\begin{align}
B_{n} := - E\left[M_{n}(\tes(X^{n});X^{n})-\tilde{E}\left\{M_{n}(\tes(X^{n});\tilde{X}^{n})\right\}\right].
\nonumber
\end{align}
Arguing as in \cite{UmeShiMasNin15}, we will show that
\begin{align}
B_{n}=p^{\ast}+q^{\ast}+o(1),
\label{hm:s.aic_formula}
\end{align}
from which, with ignoring the $o(1)$-term, we arrive at the AIC type statistics (after multiplied by $2$ as usual)
\begin{align}
\mca_{n} := 2M_{n}(\tes) + 2(p^{\ast}+q^{\ast}).
\label{ys:AIC}
\end{align}
Given more than one candidate models all belonging to our basic setup, we compute the quantity $\mca_{n}$ for each model, and then pick the one giving the minimum-$\mca_{n}$ value: as was discussed in \cite{UmeShiMasNin15} (also \cite{NinKaw16}), in some instances this criterion can be used to select the tuning parameter.

We are left to showing \eqref{hm:s.aic_formula}. By the usual mean-value theorem representation, we have: for suitable $\tes'$ and $\tes''$ both tending in $P$-probability to $\tz$,
\begin{align}
B_{n} &= -E\otimes\tilde{E}\left[\left\{M_{n}(\tes;X^{n})-M_{n}(\tz;X^{n})\right\} 
- \left\{M_{n}(\tes;\tilde{X}^{n})-M_{n}(\tz;\tilde{X}^{n})\right\}\right] \nn\\
&= -E\otimes\tilde{E}\bigg[
\bigg\{\p_{\theta}M_{n}(\tz)[\tes-\tz]+\frac{1}{2}\p_{\theta}^{2}M_{n}(\tes')[(\tes-\tz)^{\otimes 2}]\bigg\} \nn\\
&{}\qquad 
- \bigg\{\p_{\theta}M_{n}(\tz;\tilde{X}^{n})[\tes-\tz]+\frac{1}{2}\p_{\theta}^{2}M_{n}(\tes'';\tilde{X}^{n})[(\tes-\tz)^{\otimes 2}]\bigg\}\bigg].
\label{hm:bias.rep_1}
\end{align}
We observe the following two points.
\begin{itemize}
\item First we look at the weak convergence limit of the quantity inside the expectation sign ``$-E\otimes\tilde{E}$'', say $G_{n}=G_{n}(X^{n},\tilde{X}^{n})$. It follows from the sparse consistency together with the stochastic expansion \eqref{hm:linear.rep_2} that, writing $\tilde{\D}_{n}(\tz)=D_{n}\p_{\theta^{\ast}}M_{n}(\tz;\tilde{X}^{n})$, 
\begin{align}
G_{n} &= \left(\D_{n}(\tz)[\hat{w}_{n}^{\ast}] + \frac{1}{2}\Gam_{0}[\hat{w}_{n}^{\ast},\hat{w}_{n}^{\ast}] + o_{p}(1)\right)
-\left(\tilde{\D}_{n}(\tz)[\hat{w}_{n}^{\ast}] + \frac{1}{2}\Gam_{0}[\hat{w}_{n}^{\ast},\hat{w}_{n}^{\ast}] + o_{p}(1)\right)
\nn\\
&=-\left(\D_{n}(\tz)-\tilde{\D}_{n}(\tz)\right)\Gam_{0}^{-1}(\D_{n}(\tz)+\Lam^{\ast}) + o_{p}(1).
\nonumber
\end{align}
\item Second, since $(|\hat{w}_{n}|^{2})_{n}$ is assumed to be uniform integrable, we can readily deduce from the expression \eqref{hm:bias.rep_1} that $(G_{n})_{n}$ is uniform integrable as well. 
\end{itemize}
Therefore, taking the independence between $\D_{n}(\tz)$ and $\tilde{\D}_{n}(\tz)$ into account we have the convergence of moments: for an appropriate probability with expectation operator $\E$ which can integrate the weak limit of $(\D_{n}(\tz),\tilde{\D}_{n}(\tz))$, say $(\D_{0},\tilde{\D}_{0})\sim N_{p^{\ast}+q^{\ast}}(0,\Gam_{0})\otimes N_{p^{\ast}+q^{\ast}}(0,\Gam_{0})$, we obtain
\begin{align}
B_{n} &= E\otimes\tilde{E}\left\{\left(\D_{n}(\tz)-\tilde{\D}_{n}(\tz)\right)\Gam_{0}^{-1}(\D_{n}(\tz)+\Lam^{\ast})\right\} + o(1)
\nn\\
&= \E\left\{\left(\D_{0}(\tz)-\tilde{\D}_{0}(\tz)\right)\Gam_{0}^{-1}(\D_{0}(\tz)+\Lam^{\ast})\right\} + o(1)
\nn\\
&= p^{\ast}+q^{\ast} + o(1),
%&= \E\left\{\D_{0}(\tz)\Gam_{0}^{-1}\D_{0}(\tz)\right\} + o(1)
\nonumber
\end{align}
hence \eqref{hm:s.aic_formula};
we note that the constant factor $\Lam^{\ast}$ has no asymptotic contribution to the non-negligible leading term of $B_{n}$.
As the results, we can conclude the following corollary which gives a sufficient condition for \eqref{hm:s.aic_formula}, verifying the AIC type information criterion \eqref{ys:AIC}.

\begin{cor}
Suppose that $M_{n}(\tz)$ is the minus likelihood function of true data-generating model,
and that Assumptions \ref{hm:A_c}, \ref{hm:A_rc}, \ref{ys:sc}, \ref{ys:LAQ_assum}, \ref{pldi:ass1} and \ref{pldi:ass7} hold with 
$\D_{n}(\tz)\cil N_{p^{\ast}+q^{\ast}}(0,\Gam_{0})$, where $\Gam_{0}$ is a constant positive-definite matrix such that $\Gam_{n}(\tz)\cip\Gam_{0}$ given in Assumption \ref{ys:LAQ_assum}. %and $K=2+\delta$ for some $\delta>0$ in Assumption \ref{pldi:ass1} and \ref{pldi:ass7}.
Then \eqref{hm:s.aic_formula} holds.
\end{cor}

%%%%%
%%%%%
\section{Regularized estimation of discretely observed ergodic diffusion process}\label{hm:sec_ergo.diff}

In this section, we apply the results in Section \ref{hm:sec_s.asymp} to the parametric ergodic-diffusion model, keeping the basic notation of Section \ref{hm:sec_setup} in use. Let $X=(X_{t})_{t\in\mbbrp}$ be a $d$-dimensional solution process to the stochastic differential equation
\begin{align}
dX_{t}=a(X_{t},\alpha)dw_{t}+b(X_{t},\beta)dt, \quad X_{0}=x_{0},
\nonumber
\end{align}
where $\alpha=(\alpha_{1},\ldots,\alpha_{p})\in\Theta_{\alpha}\subset\mathbb{R}^{p}$, $\beta=(\beta_{1},\ldots,\beta_{q})\in\Theta_{\beta}\subset\mathbb{R}^{q}$, $a:\mathbb{R}^{d}\times\Theta_{\alpha}\rightarrow\mathbb{R}^{d}\times\mathbb{R}^{m}$, $b:\mathbb{R}^{d}\times\Theta_{\beta}\rightarrow\mathbb{R}^{d}$ and $w_{t}$, $t\in[0,T]$ is a standard Wiener process in $\mathbb{R}^{m}$ independent of the initial variable $x_{0}$.
We observe a discrete-time data $\mbX_{n}=(X_{t_{j}})_{0\leq j\leq n}$, where $t_{j}=ih_{n}$. 
We focus on a sampling design such that
\begin{align}
\text{$T_{n}:=nh_{n}\rightarrow \infty$\quad and \quad $nh_{n}^{2}\rightarrow0$\quad ($n\rightarrow\infty$).}
\nonumber
\end{align}
Let $c(x,\alpha):=a(x,\alpha)^{\otimes 2}$. We set some assumptions on $X$.
\begin{assum}\label{ys:as_Lip2}\ 
\begin{enumerate}
\item There exists a constant $C$ such that for any $x,\ y\in\mathbb{R}^{d}$
\begin{align}
\sup_{\alpha\in\overline{\Theta}_{\alpha}}|a(x,\alpha)-a(y,\alpha)|+\sup_{\beta\in\overline{\Theta}_{\beta}}|b(x,\beta)-b(y,\beta)|\leq C|x-y|.
\nonumber
\end{align}
\item $\inf_{x,\alpha}|c(x,\alpha)|>0$.
\item For all $k\in\{0,1,2\}$ and $l\in\{0,1,\ldots,4\}$, there exists a constant $C=C(k,l)$ such that
\begin{align}
\sup_{(x,\theta)\in\mathbb{R}^{d}\times\overline{\Theta}}(1+|x|)^{-C}\{|\partial_{x}^{k}\partial_{\alpha}^{l}a(x,\alpha)|+|\partial_{x}^{k}\partial_{\beta}^{l}b(x,\beta)|\}<\infty. \nonumber
\end{align}
\end{enumerate}
\end{assum}

\begin{assum}\label{ys:as_inv2}\ 
\begin{enumerate}
\item The process $X_{t}$ is ergodic for $\theta=\theta_{0}$ with invariant probability measure $\pi_{0}$:
\begin{align}
\frac{1}{T_{n}}\int_{0}^{T_{n}}f(X_{t})dt\xrightarrow{p}\int f(x)\pi_{0}(dx)
\nonumber
\end{align}
for every $f\in L^{1}(\pi_{0})$.

\item \label{ys:as_Lip2-3}If $a(x,\alpha)=a(x,\alpha_{0})$ and $b(x,\beta)=b(x,\beta_{0})$ for $\pi_{0}$-a.e. $x$, then $\alpha=\alpha_{0}$ and $\beta=\beta_{0}$. 

\item For all $k>0$, $\sup_{t}E[|X_{t}|^{k}]<\infty$.
\end{enumerate}
\end{assum}
Hereafter, we assume that Assumptions \ref{ys:as_Lip2} and \ref{ys:as_inv2} hold. We note that
\begin{align}
\frac{1}{n}\sumj f(X_{t_{j-1}})\xrightarrow{p}\int f(x)\pi_{0}(dx)
\nonumber
\end{align}
for continuously differentiable $f$ with the derivative $\partial f$ of at most polynomial growth, since
\begin{align}
E\left[\left|\frac{1}{T_{n}}\int_{0}^{T_{n}}f(X_{t})dt-\frac{1}{n}\sumj f(X_{t_{j-1}})\right|\right]\lesssim \frac{1}{n}\sumj\sup_{t_{j-1}\leq s\leq t_{j}}\sqrt{E[|X_{s}-X_{t_{j-1}}|^{2}]}\rightarrow0.
\nonumber
\end{align}

\medskip

Now, let us define a regularized estimator $\hat{\theta}_{n}=(\hat{\alpha}_{n},\hat{\beta}_{n})$ as the minimizer of \eqref{hm:Hn}: $\mathbb{H}_{n}(\theta)=M_{n}(\al,\beta)+\overline{R}_{n}^{a}(\al) +\overline{R}_{n}^{b}(\beta)$, where
\begin{align}
M_{n}(\theta)=\frac{1}{2}\sumj\left\{\log\det(c_{j-1}(\alpha))+\frac{1}{h_{n}}c_{j-1}^{-1}(\alpha)\left[(\Delta_{j}X-h_{n}b_{j-1}(\beta))^{\otimes2}\right]\right\}
%\label{ys:M-diff}
\nn
\end{align}
is a negative quasi-likelihood function up to an additive constant, $\Delta_{j}X=X_{t_{j}}-X_{t_{j-1}}$, $c_{j-1}(\alpha)=c(X_{t_{j-1}},\alpha)$ and $b_{j-1}(\beta)=b(X_{j-1},\beta)$. This quasi-likelihood has been used by many researchers, e.g., \cite{GenJac93}, \cite{Kes97}, \cite{Yos92}, and \cite{Yos11}. 
Following the route in Section \ref{hm:sec_asymptotics}, we will derive the asymptotic behavior of $\hat{\theta}_{n}$ with setting $s_{n}=1/\sqrt{n}$ and $t_{n}=1/\sqrt{T_{n}}$. 

\medskip

First, we derive the consistency $\hat{\theta}_{n}\xrightarrow{p}\theta_{0}$ by using Theorem \ref{hm:thm_consistency}. 
Note that in the present setting we have the following convergence uniformly in $\theta$ under Assumptions \ref{ys:as_Lip2} and \ref{ys:as_inv2} (see \cite[Section 6]{Yos92}):
\begin{align}
\frac{1}{n}\left\{M_{n}(\alpha,\beta_{0})-M_{n}(\alpha_{0},\beta_{0})\right\}
&=\frac{1}{2n}\sumj\bigg\{\log\frac{\det(c_{j-1}(\alpha))}{\det(c_{j-1}(\alpha_{0}))}
\nn\\
&{}\qquad+\frac{1}{h_{n}}(c_{j-1}^{-1}(\alpha)-c_{j-1}^{-1}(\alpha_{0}))[(\Delta_{j}X-h_{n}b_{j-1}(\beta_{0}))^{\otimes2}]\bigg\} \nonumber\\
&\xrightarrow{p}\frac{1}{2}\int\left\{\log\frac{\det(c(x,\alpha))}{\det(c(x,\alpha_{0}))}+\mathrm{Tr}\left(c^{-1}(x,\alpha)c(x,\alpha_{0})-I_{d}\right)\right\}\pi_{0}(dx), \nn \\
\frac{1}{T_{n}}\left\{M_{n}(\alpha,\beta)-M_{n}(\alpha,\beta_{0})\right\}&=\frac{1}{T_{n}}\sumj\Big\{\frac{h_{n}}{2}c_{i-1}^{-1}(\alpha)[b_{i-1}(\beta)^{\otimes 2}-b_{j-1}(\beta_{0})^{\otimes}] \nn \\
&\quad -c_{j-1}^{-1}(\alpha)[b_{j-1}(\beta)-b_{j-1}(\beta_{0}),\Delta_{j}X]\Big\} \nonumber \\
&\xrightarrow{p}\frac{1}{2}\int c^{-1}(x,\alpha)[(b(x,\beta)-b(x,\beta_{0}))^{\otimes2}]\pi_{0}(dx).
\nn
\end{align}
We set 
\begin{align}
\overline{M}_{0}^{a}(\al)&=\frac{1}{2}\int\bigg(\log\frac{\det(c(x,\alpha))}{\det(c(x,\alpha_{0}))}+\mathrm{Tr}\left(c^{-1}(x,\alpha)c(x,\alpha_{0})-I_{d}\right)\bigg)\pi_{0}(dx), 
\label{ys:over-Ma} \\
\overline{M}_{0}^{b}(\theta)&=\frac{1}{2}\int c^{-1}(x,\alpha)[(b(x,\beta)-b(x,\beta_{0}))^{\otimes2}]\pi_{0}(dx). 
\label{ys:over-Mb}
\end{align}
Then, we have $\overline{M}_{0}^{a}(\al)\ge 0$ and $\overline{M}_{0}^{b}(\theta)\ge 0$ with the inequalities holding if and only if $\theta=\tz$, so that the following corollary is trivial from Theorem \ref{hm:thm_consistency}.

\begin{cor}
We have $\hat{\theta}_{n}\xrightarrow{p}\theta_{0}$ under Assumptions \ref{hm:A_c}.\ref{ys:consistency-3}, \ref{ys:as_Lip2} and \ref{ys:as_inv2}.
\end{cor}

\medskip

Next, in order to show the rates of convergence $(\sqrt{n}(\hat{\alpha}_{n}-\alpha_{0}),\sqrt{T_{n}}(\hat{\beta}_{n}-\beta_{0}))=O_{p}(1)$ we apply Theorem \ref{hm:thm_rate} with $s_{n}=1/\sqrt{n}$ and $t_{n}=1/\sqrt{T_{n}}$. 
We assume Assumption \ref{hm:A_rc}.\ref{ys:R_cond.3}, and will check Assumption \ref{hm:A_rc}.\ref{ys:R_cond.1}. Conditions (a) and (c) are derived from \cite[Lemma 3, 6, 8 and 9]{Yos11}. Further, it is easy to show the conditions (b) and (d) since we have
\begin{align}
\sup_{\alpha}\left|\frac{1}{n\sqrt{h_{n}}}\p_{\al}\p_{\beta}M_{n}(\al,\beta_{0})\right|&=O_{p}(\sqrt{h_{n}}),
\nonumber \\
%\frac{1}{\sqrt{T_{n}}}\partial_{\beta}M_{n}(\aes,\beta_{0})&=-\frac{1}{\sqrt{T_{n}}}\sumj \left\{\p_{\beta}b_{j-1}^{\top}(\beta_{0})c_{j-1}^{-1}(\aes)(\Delta_{j}X-h_{n}b_{j-1}(\beta_{0}))\right\}{\color{red}=O_{p}(1)} \nn \\
\frac{1}{n}\partial_{\alpha}^{2}M_{n}(\theta_{0})&=\frac{1}{2n}\sumj\bigg\{\partial_{\alpha}^{2}\log\det(c_{j-1}(\al_{0}))
\nn\\
&{}\qquad +\frac{1}{h_{n}}\partial_{\alpha}^{2}c_{j-1}^{-1}(\al_{0})\left[(\Delta_{j}X-h_{n}b_{j-1}(\beta_{0}))^{\otimes2}\right]\bigg\} \nn \\
%&=\frac{1}{2n}\sumj\left\{\partial_{\alpha}^{2}\log|c_{j-1}(\al_{0})|+\frac{1}{h_{n}}\partial_{\alpha}^{2}c_{j-1}^{-1}(\al_{0})\left[\{\Delta_{j}X-h_{n}b_{j-1}(\beta_{0})+h_{n}(b_{j-1}(\beta_{0})-b_{j-1}(\bes))\}^{\otimes2}\right]\right\} \nn \\
%&=\frac{1}{2n}\sumj\left\{\partial_{\alpha}^{2}\log|c_{j-1}(\al_{0})|+\frac{1}{h_{n}}\partial_{\alpha}^{2}c_{j-1}^{-1}(\al_{0})\left[(\Delta_{j}X-h_{n}b_{j-1}(\beta_{0}))^{\otimes2}\right]\right\} \nn \\
%&\quad +\frac{h_{n}}{2n}\sumj\partial_{\alpha}^{2}c_{j-1}^{-1}(\al_{0})\left[(b_{j-1}(\beta_{0})-b_{j-1}(\bes))^{\otimes2}\right] \nn \\
%&\quad +\frac{1}{n}\sumj\partial_{\alpha}^{2}c_{j-1}^{-1}(\al_{0})[\Delta_{j}X-h_{n}b_{j-1}(\beta_{0}),b_{j-1}(\beta_{0})-b_{j-1}(\bes)] \nn \\
&\xrightarrow{p} \frac{1}{2}\int\left\{\partial_{\alpha}^{2}\log\det(c(x,\al_{0}))+\partial_{\alpha}^{2}c^{-1}(x,\al_{0})[c(x,\al_{0})]\right\}\pi_{0}(dx)
\nn\\
&=\bigg[ \frac{1}{2}\int {\rm Tr}\{c^{-1}(\p_{\alpha_{i}}c)c^{-1}(\p_{\al_{j}}c)(x,\al_{0})\}\pi_{0}(dx) \bigg]_{i,j=1}^{p}
=:\Gamma_{0}^{\alpha}(\al_{0}), \nn \\
\frac{1}{T_{n}}\partial_{\beta}^{2}M_{n}(\theta_{0})&=\frac{1}{2T_{n}}\sumj\left\{\frac{1}{h_{n}}c_{j-1}^{-1}(\theta_{0})\left[\p_{\beta}^{2}(\Delta_{j}X-h_{n}b_{j-1}(\beta_{0}))^{\otimes2}\right]\right\} \nn \\
%&=\frac{1}{n}\sumj c_{j-1}^{-1}(\aes)\left[(\p_{\beta}b_{j-1}(\beta_{0}))^{\otimes2}\right] \nn \\
%&\quad -\frac{1}{nh_{n}}\sumj c_{j-1}^{-1}(\aes)\left[\p_{\beta}^{2}b_{j-1}(\beta_{0})[\Delta_{j}X-h_{n}b_{j-1}(\beta_{0})]\right] \nn \\
&\xrightarrow{p}\int c^{-1}(x,\al_{0})[(\p_{\beta}b_{j-1}(\beta_{0}))^{\otimes2}]\pi_{0}(dx)=:\Gamma^{\beta}_{0}(\theta_{0}),
\nn
\end{align}
with $\lambda_{\min}(\Gamma^{\al}_{0}(\al_{0}))\wedge\lambda_{\min}(\Gamma^{\beta}_{0}(\theta_{0}))>0$;
both $\p_{\al}c(x,\al)$ and $\p_{\beta}b(x,\beta)$ are not identically zero under the assumptions.
As a result, we can conclude the following corollary:
\begin{cor}\label{ys:cor-roc}
We have $(\sqrt{n}(\hat{\alpha}_{n}-\alpha_{0}),\sqrt{T_{n}}(\hat{\beta}_{n}-\beta_{0}))=O_{p}(1)$ under Assumptions \ref{hm:A_c}.\ref{ys:consistency-3}, \ref{hm:A_rc}.\ref{ys:R_cond.3}, \ref{ys:as_Lip2} and \ref{ys:as_inv2}.
\end{cor}

%\begin{align}
%\frac{1}{\sqrt{n}}\partial_{\alpha}M_{n}(\alpha_{0},\tilde{\beta})&=\frac{1}{2\sqrt{n}}\sumj\left\{\partial_{\alpha}\log|c_{j-1}(\alpha_{0})|+\frac{1}{h_{n}}\partial_{\alpha}c_{j-1}^{-1}(\alpha_{0})\left[(\Delta_{j}X-h_{n}b_{j-1}(\tilde{\beta}))^{\otimes2}\right]\right\}; \nn \\
%\frac{1}{\sqrt{n}}\p_{\alpha}\p_{\beta}M_{n}(\al_{0},\tilde{\beta})&=-\frac{1}{\sqrt{n}}\sumj \partial_{\alpha}c_{j-1}^{-1}(\alpha_{0})\big[\p_{\beta}b_{j-1}(\tilde{\beta}),\Delta_{j}X-h_{n}b_{j-1}(\tilde{\beta})\big] \nn \\
%&=-\sqrt{h_{n}}\frac{1}{\sqrt{nh_{n}}}\sumj \partial_{\alpha}c_{j-1}^{-1}(\alpha_{0})\big[\p_{\beta}b_{j-1}(\tilde{\beta}),\Delta_{j}X-h_{n}b_{j-1}(\beta_{0})\big] \nn \\
%&\quad -(\sqrt{n}h_{n})\frac{1}{nh_{n}}\sumj \partial_{\alpha}c_{j-1}^{-1}(\alpha_{0})\big[\p_{\beta}b_{j-1}(\tilde{\beta}),h_{n}b_{j-1}(\beta_{0})-h_{n}b_{j-1}(\tilde{\beta})\big]. \nn
%\end{align}

\medskip

We proceed to the sparse consistency $P(\hat{\theta}_{n}^{\circ}=0)\rightarrow 1$; recall that we presupposed that both $\al$ and $\beta$ do have zero components. 
Let Assumptions \ref{hm:A_c}.\ref{ys:consistency-3} and \ref{hm:A_rc}.\ref{ys:R_cond.3} hold; then $(\sqrt{n}(\hat{\alpha}_{n}-\alpha_{0}),\sqrt{T_{n}}(\hat{\beta}_{n}-\beta_{0}))=O_{p}(1)$ from Corollary \ref{ys:cor-roc}, from which we deduce \eqref{ys:a-zero_tight} and \eqref{ys:b-zero_tight} with $s_{n}=1/\sqrt{n}$ and $t_{n}=/1\sqrt{T_{n}}$. We will check Assumption \ref{ys:sc}. We get
\begin{align}
&\frac{1}{n}\partial_{\alpha^{\circ}_{k'}}M_{n}(0,\ldots,0,\hat{\al}^{\circ}_{n,k'},\ldots,\hat{\al}_{n,p^{\circ}}^{\circ},\hat{\al}^{\ast}_{n},\bes)\hat{\al}^{\circ}_{n,k'}
+\frac{1}{n}\lambda_{n,k'}^{a\circ}R^{a}(\hat{\alpha}_{n,k'}^{\circ}) \nn \\
&=\bigg\{\frac{1}{\sqrt{n}}\partial_{\alpha^{\circ}_{k'}}M_{n}(0,\ldots,0,\hat{\al}^{\circ}_{n,k'},\ldots,\hat{\al}_{n,p^{\circ}}^{\circ},\hat{\al}^{\ast}_{n},\bes)\frac{\hat{\al}_{n,k'}^{\circ}}{|\hat{\al}_{n,k'}^{\circ}|}\frac{1}{|\sqrt{n}\hat{\al}_{n,k'}^{\circ}|}
\nn\\
&{}\qquad\qquad 
+\frac{1}{n}|\hat{\al}_{n,k'}^{\circ}|^{-2}\lambda_{n,k'}^{a\circ}R^{a}(\hat{\alpha}_{n,k'}^{\circ})\bigg\}|\hat{\al}^{\circ}_{n,k'}|^{2} \nn \\
&=\bigg\{O_{p}(1)\frac{1}{|\sqrt{n}\hat{\al}_{n,k'}^{\circ}|}+\frac{1}{|\sqrt{n}\hat{\al}_{n,k'}^{\circ}|^{2}}\lambda_{n,k'}^{a\circ}R^{a}(\hat{\alpha}_{n,k'}^{\circ})\bigg\}|\hat{\al}_{n,k'}^{\circ}|^{2} \nn \\
%&=\left\{\frac{1}{|s_{n}^{-1}\hat{\al}_{n,k'}^{\circ}|^{1+\gamma}}\left(O_{p}(1)+|s_{n}^{-1}\hat{\al}_{n,k'}^{\circ}|^{\gamma-1}\lambda_{n,k'}^{a\circ}|\hat{\al}^{\circ}_{n,k'}|^{\gamma}\right)\right\}|\hat{\al}_{n,k'}^{\circ}|^{2}, \nn \\
&=\bigg\{\frac{1}{|\sqrt{n}\hat{\al}_{n,k'}^{\circ}|}\left(O_{p}(1)+|\sqrt{n}\hat{\al}_{n,k'}^{\circ}|^{-1}\lambda_{n,k'}^{a\circ}R^{a}(\hat{\alpha}_{n,k'}^{\circ})\right)\bigg\}|\hat{\al}_{n,k'}^{\circ}|^{2}.
\nn 
\end{align}
Therefore, if we assume
\begin{align}
P\left(|\sqrt{n}\hat{\al}_{n,k'}^{\circ}|^{-1}\lambda_{n,k'}^{a\circ}R^{a}(\hat{\alpha}_{n,k'}^{\circ})>K\right)\rightarrow1 
\label{ys:sc-diff-a3}
\end{align}
and similarly
\begin{align}
P\left(|\sqrt{T_{n}}\hat{\beta}_{n,l'}^{\circ}|^{-1}\lambda_{n,l'}^{b\circ}R^{b}(\hat{\beta}_{n,l'}^{\circ})>K\right)\rightarrow1 
\label{ys:sc-diff-b3}
\end{align}
for any $K>0$, Assumption \ref{ys:sc} holds;
%although these conditions are not the only possibilities, they are conveniently compact.
note that the conditions \eqref{ys:sc-diff-a3} and \eqref{ys:sc-diff-b3} put limitations on specific form of the functions $R^{a}$ and $R^{b}$, and also on the stochastic orders of $\lambda_{n,k'}^{a\circ}$ and $\lambda_{n,l'}^{b\circ}$. Theorem \ref{ys:thm_sc} concludes the following.

\begin{cor}
We have $P(\hat{\theta}_{n}^{\circ}=0)\rightarrow 1$ under Assumptions \ref{hm:A_c}.\ref{ys:consistency-3}, \ref{hm:A_rc}.\ref{ys:R_cond.3}, \ref{ys:as_Lip2}, \ref{ys:as_inv2}, \eqref{ys:sc-diff-a3} and \eqref{ys:sc-diff-b3}.
\end{cor}

%we have
%\begin{align}
%\Gamma^{\al^{\circ}}_{0}(\alpha_{0})_{ij}&=\frac{1}{2}\int {\rm tr}\{c^{-1}(\p_{\alpha^{\circ}_{i}}c)c^{-1}(\p_{\al^{\circ}_{j}}c)(x,\alpha_{0})\}\pi_{0}(dx), \nn \\
%\Gamma^{\beta^{\circ}}_{0}(\theta_{0})&=\int\left\{ c^{-1}(x,\al_{0})[(\p_{\beta^{\circ}}b_{j-1}(\beta_{0}))^{\otimes2}]\right\}\pi_{0}(dx). \nn
%\end{align}

%\begin{align}
%&\frac{1}{n}\partial_{\alpha^{\circ}_{k'}}M_{n}(0,\ldots,0,\hat{\alpha}^{\circ}_{n,k'},\ldots,\hat{\alpha}^{\circ}_{n,p^{\circ}},\aes,\bes) \nn \\
%&=\frac{1}{2n}\sumj\left\{\partial_{\alpha^{\circ}_{k'}}\log|c_{j-1}(0,\ldots,0,\hat{\alpha}_{n,k'}^{\circ},\ldots,\hat{\alpha}_{n,p^{\circ}},\aes^{\ast})|\right. \nn \\
%&\quad\left.+\frac{1}{h_{n}}\partial_{\alpha^{\circ}_{k'}}c_{j-1}^{-1}(0,\ldots,0,\hat{\alpha}_{n,k'}^{\circ},\ldots,\hat{\alpha}_{n,p^{\circ}},\aes^{\ast})\left[(\Delta_{j}X-h_{n}b_{j-1}(\bes))^{\otimes2}\right]\right\} \nn 
%\end{align}
%\begin{align}
%&\frac{1}{nh_{n}}\partial_{\beta^{\circ}_{l'}}M_{n}(\aes,0,\ldots,0,\hat{\beta}_{n,l'}^{\circ},\ldots,\hat{\beta}_{n,q^{\circ}},\bes^{\ast}) \nn \\
%&=\frac{-1}{nh_{n}}\sumj c_{j-1}(\aes)[\Delta_{j}X-h_{n}b_{j-1}(\aes,0,\ldots,0,\hat{\beta}_{n,l'}^{\circ},\ldots,\hat{\beta}_{n,q^{\circ}},\bes^{\ast}),\p_{\beta^{\circ}_{l'}}b_{j-1}(\aes,0,\ldots,0,\hat{\beta}_{n,l'}^{\circ},\ldots,\hat{\beta}_{n,q^{\circ}},\bes^{\ast})] \nn 
%\end{align}

\medskip

We now apply Corollary \ref{ys-cor_AL} to derive the asymptotic non-degenerate distribution of
\begin{align}
(\hat{u}_{n}^{\ast},\hat{v}_{n}^{\ast})=\left(\sqrt{n}(\hat{\alpha}_{n}^{\ast}-\alpha_{0}^{\ast}),\sqrt{T_{n}}(\hat{\beta}_{n}^{\ast}-\beta_{0}^{\ast})\right).
\nonumber
\end{align} 
We use the same notation as in Section \ref{hm:sec_ad.nz}; recall that
\begin{align}
\mbbm_{n}^{0}(u^{\ast},v^{\ast}) :=M_{n}(\aes^{\circ},\al_{0}^{\ast}+s_{n}u^{\ast},\bes^{\circ},\beta_{0}^{\ast}+t_{n}v^{\ast})-M_{n}(\aes^{\circ},\al^{\ast}_{0},\bes^{\circ},\beta^{\ast}_{0}).
\nonumber
\end{align}
By means of the Taylor expansion around $(u^{\ast},v^{\ast})=(0,0)$ together with the property $P(\hat{\theta}_{n}^{\circ}=0)\rightarrow 1$ and the martingale central limit theorem, we obtain
\begin{align}
\mathbb{M}_{n}^{0}(u^{\ast},v^{\ast}) 
%&=\frac{1}{2}\sumj\left\{\log\det\left(c_{j-1}\left(0,\al_{0}^{\ast}+\frac{u^{\ast}}{\sqrt{n}}\right)\right)+\frac{1}{h_{n}}c_{j-1}^{-1}\left(0,\al_{0}^{\ast}+\frac{u^{\ast}}{\sqrt{n}}\right)\left[\left(\Delta_{j}X-h_{n}b_{j-1}\left(0,\beta_{0}^{\ast}+\frac{v^{\ast}}{\sqrt{nh_{n}}}\right)\right)^{\otimes2}\right]\right\} \nn \\
%&\quad-\frac{1}{2}\sumj\left\{\log\det(c_{j-1}(\alpha_{0}))+\frac{1}{h_{n}}c_{j-1}^{-1}(\alpha_{0})\left[(\Delta_{j}X-h_{n}b_{j-1}(\beta_{0}))^{\otimes2}\right]\right\} \nn \\
%&=\frac{1}{2\sqrt{n}}\sumj\left\{\p_{\al^{\ast}}\log\det(c_{j-1}(\al_{0}))+\frac{1}{h_{n}}\p_{\al^{\ast}}c_{j-1}^{-1}(\al_{0})\left[\left(\Delta_{j}X-h_{n}b_{j-1}\left(\beta_{0}\right)\right)^{\otimes2}\right]\right\}[u^{\ast}] \nn \\
%&\quad+\frac{1}{2\sqrt{nh_{n}}}\sumj\left(\frac{1}{h_{n}}c_{j-1}^{-1}\left(\al_{0}\right)\left[\p_{\beta^{\ast}}\left(\Delta_{j}X-h_{n}b_{j-1}\left(\beta_{0}\right)\right)^{\otimes2}\right]\right)[v^{\ast}]+\frac{1}{4n}\sumj\Big\{\p_{\al^{\ast}}^{2}\log\det(c_{j-1}(\al_{0}))[u^{\ast},u^{\ast}] \nn \\
%&\quad+\frac{1}{h_{n}}\p_{\al^{\ast}}^{2}c_{j-1}^{-1}(\al_{0})[{u^{\ast}}^{\otimes2},\left(\Delta_{j}X-h_{n}b_{j-1}(\beta_{0})\right)^{\otimes2}]\Big\} \nn \\
%&\quad+\frac{1}{4nh_{n}}\sumj\left(\frac{1}{h_{n}}c_{j-1}^{-1}\left(\al_{0}\right)\left[\p_{\beta^{\ast}}^{2}\left(\Delta_{j}X-h_{n}b_{j-1}\left(\beta_{0}\right)\right)^{\otimes2}\right]\right)[v^{\ast},v^{\ast}]+o_{p}(1) \nn \\
&=\frac{1}{2\sqrt{n}}\sumj\bigg\{\p_{\al^{\ast}}\log\det(c_{j-1}(\al_{0}))
\nn\\
&{}\qquad
+\frac{1}{h_{n}}\p_{\al^{\ast}}c_{j-1}^{-1}(\al_{0})\left[\left(\Delta_{j}X-h_{n}b_{j-1}\left(\beta_{0}\right)\right)^{\otimes2}\right]\bigg\}[u^{\ast}] \nn \\
&{}\quad+\frac{1}{2\sqrt{nh_{n}}}\sumj\left(\frac{1}{h_{n}}c_{j-1}^{-1}\left(\al_{0}\right)\left[\p_{\beta^{\ast}}\left(\Delta_{j}X-h_{n}b_{j-1}\left(\beta_{0}\right)\right)^{\otimes2}\right]\right)[v^{\ast}]
\nn\\
&{}\quad +\frac{1}{4}\int\Big\{\p_{\al^{\ast}}^{2}\log\det(c(x,\al_{0}))[u^{\ast},u^{\ast}] +\p_{\al^{\ast}}^{2}c^{-1}(x,\al_{0})[{u^{\ast}}^{\otimes2},c(x,\al_{0})]\Big\}\pi_{0}(dx) \nn\\
&{}\quad+\frac{1}{2}\int c^{-1}(x,\al_{0})[\p_{\beta^{\ast}}b(x,\beta_{0})[v^{\ast}],\p_{\beta^{\ast}}b(x,\beta_{0})[v^{\ast}]]\pi_{0}(dx)+o_{p}(1) \nn \\
%&\quad+\frac{1}{2nh_{n}}\sumj\left(\frac{1}{h_{n}}c_{j-1}^{-1}\left(\al_{0}\right)\left[\p_{\beta^{\ast}}^{2}\left(\Delta_{j}X-h_{n}b_{j-1}\left(\beta_{0}\right)\right)^{\otimes2}\right]\right)[v^{\ast},v^{\ast}]+o_{p}(1) \nn \\
%&+\frac{1}{2n}\sumj\int_{0}^{1}\p_{\al^{\ast}}^{2}\log\left|c_{j-1}\left(0,\al_{0}^{\ast}+\frac{u^{\ast}}{\sqrt{n}}\xi\right)\right|d\xi\left[u^{\ast},u^{\ast}\right]+\frac{1}{2n\sqrt{h}}\sumj\int_{0}^{1}\frac{1}{h_{n}}\p_{\theta^{\ast}}^{2}
&\xrightarrow{\mathcal{L}}\Delta_{0}(\theta_{0})[w^{\ast}]+\frac{1}{2}\Gamma_{0}[w^{\ast},w^{\ast}],
\nn
\end{align}
where
\begin{align}
\Delta_{0}(\theta)&\sim N\left(0, \diag(\Delta_{0}^{11}(\theta_{0}),\Delta_{0}^{22}(\theta_{0}))\right),\nn \\
\Delta_{0}^{11}(\theta_{0})&=\frac{1}{2}\int\left\{\p_{\al^{\ast}}^{2}\log\det(c(x,\al_{0}))+\p_{\al^{\ast}}^{2}c^{-1}(x,\al_{0})[c(x,\al_{0})]\right\}\pi_{0}(dx), \nn \\
\Delta_{0}^{22}(\theta_{0})&=\int c^{-1}(x,\al_{0})[\p_{\beta^{\ast}}b(x,\beta),\p_{\beta^{\ast}}b(x,\beta_{0})]\pi_{0}(dx), \nn \\
\Gamma_{0}&=\diag\left(\frac{1}{2}\int\left\{\p_{\al^{\ast}}^{2}\log\det(c(x,\al_{0}))+\p_{\al^{\ast}}^{2}c^{-1}(x,\al_{0})[c(x,\al_{0})]\right\}\pi_{0}(dx),\right. \nn \\
&\quad\qquad\left.\int c^{-1}(x,\al_{0})[\p_{\beta^{\ast}}b(x,\beta_{0}),\p_{\beta^{\ast}}b(x,\beta_{0})]\pi_{0}(dx)\right). \nn
\end{align}
Hence, Assumptions \ref{ys:LAQ_assum}.\ref{ys:and-LAQ1} and \ref{ys:LAQ_assum}.\ref{ys:and-LAQ2} hold. Obviously, Assumption \ref{hm:A_rc}.\ref{ys:R_cond.3} is implied by Assumption \ref{ys:LAQ_assum}.\ref{ys:and-LAQ3}. The following corollary is then trivial from Corollary \ref{ys-cor_AL}:
\begin{cor}
Grant Assumptions \ref{hm:A_c}.\ref{ys:consistency-3}, \ref{ys:LAQ_assum}.\ref{ys:and-LAQ3} with $\Lam^{\ast}=(\Lam^{a \ast},\Lam^{b \ast})$ given by \eqref{hm:Lam.a_def} and \eqref{hm:Lam.b_def} being nonrandom, \ref{ys:as_Lip2}, \ref{ys:as_inv2}, \eqref{ys:sc-diff-a3}, and \eqref{ys:sc-diff-b3}.
Then, the asymptotic distribution of $(\hat{u}_{n}^{\ast},\hat{v}_{n}^{\ast})$ is
\begin{align}
(\hat{u}_{n}^{\ast},\hat{v}_{n}^{\ast})\xrightarrow{\mathcal{L}}N_{p^{\ast}+q^{\ast}}(-\Gamma_{0}^{-1}\Lambda^{\ast},\Gamma_{0}^{-1}).
\nn
\end{align}
\end{cor}

\medskip

Finally, we look at the tail-probability estimate of $(\hat{u}_{n},\hat{v}_{n})$ by means of Theorem \ref{thm:pldi}. 
For this purpose, we strengthen the ergodicity assumption and Assumption \ref{ys:as_inv2}.\ref{ys:as_Lip2-3}:
\begin{assum} \label{ys:non-deg}
\begin{enumerate}
\item 
$X$ is exponentially strong-mixing: there exists a constant $C>0$ such that
\begin{align}
\sup_{t\in\mbbrp}\sup_{A\in\sig(X_{r}: r\le t) \atop B\in\sig(X_{r}: r\ge t+h)} \left| P(A\cap B) - P(A)P(B) \right| \le C^{-1}e^{-Ch},\qquad h>0.
\nonumber
\end{align}

\item There exist positive constants $\chi^{a}$ and $\chi^{b}$, such that
\begin{enumerate}
\item $\overline{M}_{0}^{a}(\al)\geq \chi^{a}|\al-\al_{0}|^{2}$ for all $\al\in\Theta_{\al}$,
\item $\overline{M}_{0}^{b}(\al_{0},\beta)\geq \chi^{b}|\beta-\beta_{0}|^{2}$ for all $\beta\in\Theta_{\beta}$,
\end{enumerate}
where $\overline{M}_{0}^{a}(\al)$ and $\overline{M}_{0}^{b}(\theta)$ are given by \eqref{ys:over-Ma} and \eqref{ys:over-Mb}.
\end{enumerate}
\end{assum}

Thanks to Assumption \ref{ys:non-deg} and the results in \cite[Section 6]{Yos11}, we see that all the conditions in Assumption \ref{pldi:ass1} hold. Also noted is that Assumption \ref{hm:A_c}.\ref{ys:consistency-3} is implied by Assumption \ref{pldi:ass7}. Hence we obtain the following result:

\begin{cor}
For any $L>0$, \eqref{ys:joint-PLDI} holds under Assumptions \ref{pldi:ass7}, \ref{ys:as_Lip2}, \ref{ys:as_inv2} and \ref{ys:non-deg}. 
Additionally if we assume Assumption \ref{ys:LAQ_assum}.\ref{ys:and-LAQ3} with $\Lam^{\ast}=(\Lam^{a \ast},\Lam^{b \ast})$ being nonrandom, \eqref{ys:sc-diff-a3}, and \eqref{ys:sc-diff-b3}, then \eqref{eq:p.g.} holds for $\hat{u}_{0}^{\circ}=0$, $\hat{v}_{0}^{\circ}=0$, $(\hat{u}_{0}^{\ast},\hat{v}_{0}^{\ast})\sim N_{p^{\ast}+q^{\ast}}(-\Gamma_{0}^{-1}\Lambda^{\ast},\Gamma_{0}^{-1})$, and all continuous $f:\mbbr^{p+q}\to\mbbr$ of at most polynomial growth.
\end{cor}

\medskip

We thus conclude that, for example with a bridge type regularization terms for $\gam<1$ it is possible to make the asymptotic variances of the zero parameters degenerate with leaving that of non-zero part unbiased in the limit: see the discussions in Section \ref{hm:sec_tp.selection} together with the conditions \eqref{ys:sc-diff-a3} and \eqref{ys:sc-diff-b3}.

\medskip

Let us again note that non-sparse asymptotic is also in our scope (Section \ref{hm:sec_st.asymp.}). For example, as a formal variant of the drift estimation of a diffusion process we may consider the ridge-type estimation of the drift coefficient with leaving the diffusion one completely unknown. That is, for the model
\begin{align}
dX_{t}=\sum_{l=1}^{q}\beta_{l}b_{l}(X_{t})dt + a(X_{t})dw_{t},
\nonumber
\end{align}
we could consider the regularized version of the least-squares type contrast function (e.g. \cite{Mas05jjss} and the references therein):
\begin{align}
\mbbh_{n}(\beta) := \frac{1}{2h_{n}}\sumj \bigg(\D_{j}X - h_{n}\sum_{l=1}^{q}\beta_{l}b_{l}(X_{t_{j-1}})\bigg)^{2} + \sum_{l=1}^{q}\lam_{n,l}\beta_{l}^{2}.
\nonumber
\end{align}
The associated regularized estimator is given by
\begin{align}
\bes = \bigg( \sumj b_{j-1}^{\otimes 2} + 2\diag(\lam_{n,1},\dots,\lam_{n,q}) \bigg)^{-1}\sumj (\D_{j}X)b_{j-1}
\nonumber
\end{align}
with $\lam_{n,k}>0$ a.s., where $b_{j-1}=(b_{1}(X_{t_{j-1}}),\dots,b_{q}(X_{t_{j-1}}))$, would be more variance-stabilized if the computation of the random matrix $(\sumj b_{j-1}^{\otimes 2})^{-1}$ is unstable.

%%%%%
%%%%%
\appendix
\section{Polynomial type tail-supremum probability estimate}\label{appendix-A}

For convenience of reference, we here state a version of the PLDI of \cite[Theorem 1]{Yos11}. 
We need to introduce some notation. Given sets $\mathcal{T}$ and $K\subset \Theta \times \mathcal{T}$, 
we let $\xi_{0}:=(\theta_{0},\tau_{0})\in K$ be the true value of parameter $\xi$. 
Let $\ep_{n}=\ep_{n}(\theta_{0})\to 0$ be a nonrandom positive sequence. 
For contrast functions $\mathbb{H}_{n}:\Omega\times\Theta\times\mathcal{T}\to\mathbb{R}$, we define random functions
\begin{align}
\mathbb{M}_{n}(u,\tau;\theta_{0})&:=\mathbb{H}_{n}\bigl(\theta_{0}+\ep_{n}u,\tau\bigr)-\mathbb{H}_{n}(\theta_{0},\tau), \nn\\
\mathbb{Y}_{n}(\theta,\tau;\xi_{0})&:=-\ep_{n}^{2}\bigl\{\mathbb{H}_{n}(\theta,\tau)-\mathbb{H}_{n}(\theta_{0},\tau)\bigr\}.
\nonumber
\end{align}
Also, let $(\theta,\tau)\mapsto\mathbb{Y}_{0}(\theta,\tau;\xi_{0})$ be a random function. 
We consider the LAQ representation of $\mathbb{M}_{n}$:
\begin{align}
\mathbb{M}_{n}(u,\tau;\theta_{0})=\Delta_{n}(\tau;\xi_{0})[u]-\frac{1}{2}\Gamma_{0}(\tau;\xi_{0})[u,u]+r_{n}(u,\tau;\xi_{0})
\nn
%\label{plaq_rep}
\end{align}
for $(u,\tau)\in \{u\in \mathbb{R}^{p}:\ \theta_{0}+\ep_{n}u\in \Theta\}\times \mathcal{T}$, where 
$\Delta_{n}(\tau;\xi_{0})\in \mathbb{R}^{p}$, $\Gamma_{0}(\tau;\xi_{0})\in \mathbb{R}^{p}\otimes \mathbb{R}^{p}$ 
and $r_{n}(u,\tau;\xi_{0})\in \mathbb{R}$ are random variables. 
Finally, let $\zeta\in(0,1)$ and $U_{n}(r,\xi_{0}):=\{u \in \mathbb{R}^{p}:r\leq \vert u\vert \leq \ep_{n}^{-(1-\zeta)}\}$. 
We now introduce some conditions.
\begin{itemize}
\item[{\rm [A1]}] $\exists \rho_{1}>0,\ \forall L>0,$
\begin{align}
\sup_{r>0}\sup_{n>0}r^{L}P\bigg(\sup_{(u,\tau)\in U_{n}(r,\xi_{0})\times \mathcal{T}}\frac{\vert r_{n}(u,\tau;\xi_{0})\vert}{1+\vert u\vert^{2}}\geq r^{-\rho_{1}}\bigg)<\infty. \nonumber
\end{align}

\item[{\rm [A2]}] $\forall L>0,$
\begin{align}
\sup_{r>0}r^{L}P\bigg(
\forall(u,\tau)\in \mathbb{R}^{p}\times \mathcal{T},\ \ 
r^{-\rho_{1}} \vert u\vert^{2}> \frac{\Gamma_{0}(\tau;\xi_{0})}{4}[u,u]
\bigg)<\infty. \nonumber
\end{align}

\item[{\rm [A3]}] $\forall \xi_{0}\in K,\ \exists \chi(\xi_{0})>0:{\rm random\ variable},\ \exists \rho=\rho(\xi_{0})>0,\ \forall \theta\in \Theta,\ \forall \tau \in \mathcal{T},$
\begin{align}
\mathbb{Y}_{0}(\theta,\tau;\xi_{0})=\mathbb{Y}_{0}(\theta,\tau;\xi_{0})-\mathbb{Y}_{0}(\theta_{0},\tau;\xi_{0})\leq -\chi(\xi_{0})\vert \theta-\theta_{0}\vert^{\rho}. \nonumber
\end{align}

\item[{\rm [A4]}]  $\zeta\in(0,1),\ \rho_{1}\in(0,1),\ \zeta \rho<\rho_{2},\ \nu\in[0,\infty),\ 1-2\nu-\rho_{2}>0$.

\item[{\rm [A5]}] 
$\forall L>0,$
\begin{align}
\sup_{r>0}r^{L}P\left(\chi(\xi_{0})\leq r^{-(\rho_{2}-\zeta \rho)}\right)<\infty. \nonumber
\end{align}

\item[{\rm [A6]}] 
$\forall L>0,\ N_{1}:=L(1-\rho_{1})^{-1},\ N_{2}:=L(1-2\nu-\rho_{2})^{-1},$
\begin{align}
%&\qquad \qquad \qquad \sup_{\xi_{0}\in K}\sup_{n>0}E\Big[\Bigl(\sup_{\tau \in T}\vert \Delta_{n}(\tau;\xi_{0})\vert \Bigr)^{N_{1}}\Big]<\infty; \nonumber \\
&\sup_{n>0}E\bigg\{\bigg(\sup_{\tau \in T}\vert \Delta_{n}(\tau;\xi_{0})\vert \bigg)^{N_{1}}\bigg\}<\infty, \nn \\
&\sup_{n>0}E\bigg\{\bigg(\sup_{(\theta,\tau)\in\Theta\times\mathcal{T}}
\ep_{n}^{-2(1/2-\nu)}\big| \mathbb{Y}_{n}(\theta,\tau;\xi_{0})-\mathbb{Y}_{0}(\theta,\tau;\xi_{0})\bigr| \bigg)^{N_{2}}\bigg\}<\infty. \nonumber
\end{align}

\end{itemize}

The conditions [A2], [A3] and [A5] may be automatic when 
$\theta\mapsto \mathbb{Y}_{0}(\theta,\tau;\xi_{0})$ is smooth enough, and $\Gamma_{0}(\tau;\xi_{0})>0$ and $\chi(\xi_{0})>0$ are nonrandom 
uniformly in the nuisance element $\tau$. 
This is the case for many kinds of ergodic models. We also introduce:
\begin{itemize}
\item[{\rm [B1]}] $\Gamma_{0}(\tau;\xi_{0})$ is deterministic and positive-definite uniformly in $\tau\in\mathcal{T}$ and $\xi_{0}\in K$.
\item[{\rm [B2]}]  $\exists \chi=\chi(\xi_{0})>0$ : nonrandom, $\exists \rho=\rho(\xi_{0})>0,\ \forall \theta\in \Theta,\ \forall \tau\in \mathcal{T},$
\begin{align}
\mathbb{Y}_{0}(\theta,\tau;\xi_{0})\leq -\chi \vert \theta-\theta_{0}\vert^{\rho}. \nonumber
\end{align}
\end{itemize}
%Let $\hat{u}_{n}:=\ep_{n}^{-1}(\hat{\theta}_{n}-\theta_{0})$.
Then the following theorem holds:

\begin{thm}[Yoshida \cite{Yos11}, Theorems 1 and 3(a)]\label{thm_NY2011}
Assume either ``[A1]--[A6]'', or ``[A1], [A4], [A6], [B1], and [B2]''. 
Then, for every $L>0$
\begin{align}
\sup_{r>0}\sup_{n>0}r^{L}P\left(\big|\ep_{n}^{-1}(\hat{\theta}_{n}-\theta_{0})\big|\ge r\right)<\infty.
\nonumber
\end{align}
\end{thm}
The incorporation of the ``nuisance'' element $\tau$ plays an important role when the estimator has components converging at different rates 
such as the case of ergodic diffusion treated in Section \ref{hm:sec_ergo.diff}. 
Then, the PLDI would be derived by repeated use of Theorem \ref{thm_NY2011}; 
we refer to \cite[Section 5]{Yos11} and also \cite[Section 4]{Mas_LM} for a detailed account.

\bigskip

\noindent {\bf{Acknowledgements}} The authors are grateful to the referee for his helpful comments and suggestion, all of which led to substantial improvement. This work was supported by JST, CREST (H. Masuda) and JSPS (Y. Shimizu).

%%%%%
%\bibliographystyle{abbrv}
%\bibliography{YS_bibs}

%%%%%
%%%%%

\end{document}